\documentclass[leqno,12pt]{amsart}

\usepackage{amssymb, amsmath, mathrsfs}
\usepackage{amscd}              
\usepackage[all]{xy}            
\usepackage[T1]{fontenc}
\usepackage{bbold}              

\xyoption{curve}
\CompileMatrices

\newcommand{\A}{\mathbb{A}}

\newcommand{\Aff}{\text{Aff}}
\newcommand{\As}{\mathscr{A}}

\newcommand{\BGm}{\operatorname{B\Gm}}

\newcommand{\C}{\mathbb{C}}
\newcommand{\Cc}{{\mathcal{C}}}

\newcommand{\clutch}{\mathfrak{c}}

\newcommand{\comp}{\circ}

\newcommand{\dual}{\vee}

\newcommand{\exchange}{\text{exchange}}

\newcommand{\E}{{\mathcal{E}}}

\newcommand{\Ens}{\text{Ens}}

\newcommand{\ev}{\text{ev}}

\newcommand{\Gm}{\mathbb{G}_m}

\newcommand{\Hom}{\text{Hom}}
\newcommand{\Ic}{\mathcal{I}}
\newcommand{\id}{\text{id}}

\newcommand{\injto}{\hookrightarrow}

\newcommand{\isomorph}{\cong}           
\newcommand{\isomto}{\overset{\sim}{\rightarrow}}  
\newcommand{\I}{{\mathcal{I}}}

\newcommand{\Isom}{\text{Isom}}
\newcommand{\Isoms}{\text{{\it $\mathscr{I}$\!\!som}}}

\newcommand{\It}{\widetilde{I}}

\newcommand{\kb}{\mathbb{k}}

\newcommand{\Ll}{{\mathcal{L}}}                  
\newcommand{\M}{{\mathcal{M}}}
\newcommand{\Moz}{\Mt_{0,2}}
\newcommand{\Mb}{\overline{M}}

\newcommand{\mo}{{-1}}
\newcommand{\Mt}{\widetilde{\M}}
\newcommand{\N}{{\mathbb{N}}}                    

\newcommand{\Num}{\text{Num}}

\newcommand{\ob}{\text{ob}}
\newcommand{\Oo}{{\mathcal{O}}}                  
\newcommand{\Oplus}{\bigoplus}

\newcommand{\Pc}{{\mathcal{P}}}
\newcommand{\PD}{\operatorname{PD}}

\newcommand{\phit}{\tilde{\phi}}

\newcommand{\Pp}{{\mathbb{P}}}                   
\newcommand{\pr}{\text{pr}}
\newcommand{\Proj}{\text{Proj}\,}
\newcommand{\psit}{\tilde{\psi}}
\newcommand{\Q}{{\mathbb{Q}}}                    
\newcommand{\Qc}{{\mathcal{Q}}}

\newcommand{\Sb}{\mathbb{S}}

\newcommand{\Sing}{\text{Sing}}
\newcommand{\Spec}{\text{Spec}\, }               
\newcommand{\Ss}{\mathscr{S}}
\newcommand{\st}{\ \ |\ \ }
\newcommand{\sta}{\operatorname{st}}
\newcommand{\stt}{\widetilde{\sta}}
\newcommand{\Sym}{\text{Sym}\, }
\newcommand{\tensor}{\otimes}

\newcommand{\Tails}{\mathscr{T}}
\newcommand{\Tensor}{\bigotimes}

\newcommand{\Ws}{\mathscr{W}}
\newcommand{\x}{{\text{{\boldmath $x$}}}}    

\newcommand{\Xs}{{\mathscr{X}}}

\newcommand{\Ys}{{\mathscr{Y}}}

\newcommand{\Z}{{\mathbb{Z}}}                    
\newcommand{\Zar}{\text{Zar}}

\newcommand{\Zs}{\mathscr{Z}}

\newcommand{\beq}{\begin{eqnarray*}}
\newcommand{\eeq}{\end{eqnarray*}}
\newcommand{\bmul}{\begin{multline*}}
\newcommand{\emul}{\end{multline*}}
\newcommand{\broman}{\def\theenumi{\roman{enumi}}
                     \def\labelenumi{(\theenumi)}
                     \begin{enumerate}}
\newcommand{\eroman}{\end{enumerate}}
\newcommand{\balph}{\def\theenumi{\alph{enumi}}
                     \def\labelenumi{(\theenumi)}
                     \begin{enumerate}}
\newcommand{\ealph}{\end{enumerate}}
\newcommand{\barabic}{\def\theenumi{\arabic{enumi}}
                    \def\labelenumi{(\theenumi)}
                    \begin{enumerate}}
\newcommand{\earabic}{\end{enumerate}}

\newcommand{\bdefin}{\begin{definition}}
\newcommand{\edefin}{\end{definition}}
\newcommand{\bex}{\begin{example}}
\newcommand{\eex}{\end{example}}
\newcommand{\brem}{\begin{remark}}
\newcommand{\erem}{\end{remark}}
\newcommand{\bcl}{\begin{claim}}
\newcommand{\ecl}{\end{claim}}
\newcommand{\blem}{\begin{lemma}}
\newcommand{\elem}{\end{lemma}}
\newcommand{\bpr}{\begin{proposition}}
\newcommand{\epr}{\end{proposition}}
\newcommand{\bthm}{\begin{theorem}}
\newcommand{\ethm}{\end{theorem}}
\newcommand{\bco}{\begin{corollary}}
\newcommand{\eco}{\end{corollary}}
\newcommand{\bcon}{\begin{conjecture}}
\newcommand{\econ}{\end{conjecture}}
\newcommand{\bpf}{\begin{proof}}
\newcommand{\epf}{\end{proof}}
\newcommand{\bma}{\begin{pmatrix}}
\newcommand{\ema}{\end{pmatrix}}
\newcommand{\bca}{\begin{cases}}
\newcommand{\eca}{\end{cases}}

\newtheorem{theorem}{Theorem}[section]
\newtheorem{proposition}[theorem]{Proposition}
\newtheorem{lemma}[theorem]{Lemma}
\newtheorem{corollary}[theorem]{Corollary}
\newtheorem{claim}[theorem]{Claim}

\theoremstyle{definition}
\newtheorem{definition}[theorem]{Definition}
\newtheorem{remark}[theorem]{Remark}
\newtheorem{example}[theorem]{Example}
\newtheorem{conjecture}[theorem]{Conjecture}

\makeatletter
\def\overunderbraces #1#2#3{{%
 \baselineskip\z@skip \lineskip4\p@ \lineskiplimit4\p@
 \displaystyle  
 \setbox\z@\vbox{\ialign{&\hfil${}##{}$\hfil\cr
   \global\let\br\br@label #1\cr 
   \global\let\br\br@down #1\cr   
   #2\cr 
 }}
 \dimen@-\ht\z@ 
 \setbox\z@\vbox{\ialign{&\hfil${}##{}$\hfil\cr
   \global\let\br\br@label #1\cr 
   \global\let\br\br@down #1\cr   
   #2\cr 
   \global\let\br\br@up #3\cr 
   \global\let\br\br@label #3\cr   
 }}
 \advance\dimen@\ht\z@ 
 \lower\dimen@\hbox{\box\z@} 
}}
\def\br@up#1#2{\multispan{#1}\upbracefill}
\def\br@down#1#2{\multispan{#1}\downbracefill}
\def\br@label#1#2{\multispan{#1}\hidewidth $#2$\hidewidth}
\makeatother

\begin{document}

\title[A Cyclic Operad]
{A Cyclic Operad in the Category of Artin Stacks and
 Gravitational Correlators}
\author[Ivan Kausz]{Ivan Kausz}
\dedicatory{To Zsuzsanna with love}
\date{\today}
\email{ivan.kausz@mathematik.uni-regensburg.de}
\begin{abstract}
We define an Artin stack which may be considered
as a substitute for the non-existing (or empty) moduli
space of stable two-pointed curves of genus zero.
We show that this Artin stack can be viewed as the first term
of a cyclic operad in the category of stacks.
Applying the homology functor we obtain a linear
cyclic operad. 
We formulate conjectures which assert that cohomology of a smooth
projective variety has the structure of an algebra
over this homology operad and that gravitational quantum cohomology
can naturally expressed in terms of this algebra.
As a test for these conjectures we show how certain well-known relations
between gravitational correlators can be deduced from them.
\end{abstract}

\maketitle

\tableofcontents


\section{Introduction}

It is well known that some aspects of 
tree-level Gromov-Witten theory of a smooth
projective variety can be summarised by saying that its
cohomology with values in a suitable ring has the structure
of an algebra over the  operad $(A_*(\Mb_{0,n+1}))_{n\geq 2}$.
In the papers 
\cite{Manin} \S VI. 1.3, 7.6, \cite{Manin-Zograf} \S 7, \cite{LM} \S 4.1
Yu. Manin pointed out  
that it is desirable to construct an extended operad 
$(A_*(\Mt_{0,n+1}))_{n\geq 1}$
which includes the $n=1$ term, such that tree-level Gromov-Witten invariants
with gravitational descendants of a smooth projective variety
can naturally be formulated in the context of algebras over that operad.
In this paper we achieve that goal by defining stacks
$$
\Moz:=\BGm\times\Ws
\quad\text{and}\quad
\Mt_{0,n+1}:=\Mb_{0,n+1}\times\Ws^{n+1}
\quad\text{for $n\geq 2$,}
$$
where $\Ws$ is a zero-dimensional Artin-stack parameterising 
modifications of a fixed one-pointed curve (cf. \S \ref{section W}
for a precise definition of $\Ws$).

In order to be the first term of an operad, $\Moz$ has to
carry the structure of a semigroup. 
In \S \ref{section zollstock} we prove that $\Moz$ is the moduli stack
for objects which we call ``zollstocks'' (the name is derived from the German 
word for ``folding rules''). 
A zollstock is a two-pointed prestable curve of genus zero
whose dual graph is linear and whose marked points are on the two extremal
components. From this interpretation the semigroup structure of 
$\Moz$ becomes evident: Given two zollstocks, one obtains a third
zollstock by clutching together the first and the second at marked points.
In \S \ref{section operad} we show more generally that $(\Mt_{0,n+1})_{n\geq 1}$ 
has a natural structure of a cyclic operad in the category of Artin-stacks.

Passing to homology one obtains the linear cyclic operad 
$(A_*(\Mt_{0,n+1}))_{n\geq 1}$.
Let $R$ be the intersection ring of the zero-dimensional Artin stack $\Ws$.
Then we have isomorphisms
\begin{eqnarray*}
A_*(\Moz) &\isomorph& R[t_0]
\qquad\text{and}
\\
A_*(\Mt_{0,n+1}) &\isomorph&
A_*(\Mb_{0,n+1})\tensor R^{\tensor (n+1)}
\qquad\text{for $n\geq 2$.}
\end{eqnarray*}
The concrete structure of the graded ring
$R$ is given in \S\ref{section ring of W} on page \pageref{R}.

There is a natural involution on the stack of zollstocks (and
thus on $\Moz$), which
comes from interchanging the two marked points.
This involution  induces an involution $\iota$ on 
$R[t_0]$ for which we give an explicit formula
in Proposition \ref{iota}.

The semigroup structure on $\Moz$ induces a (non-commutative and
non-unital) ring structure
on $A_*(\Moz)$. Thus the group isomorphism  
$A_*(\Moz)\isomorph R[t_0]$
induces a second product on $R[t_0]$, which we denote by $\odot$.
In Theorem \ref{structure} we show that  $R[t_0]$
has a basis consisting of elements
$$
t_0^{k_1}\odot\dots\odot t_0^{k_r}
$$
where $r\geq 1$ and $k_i\geq 0$.

Let $V$ be a smooth projective variety.
For $g,n\geq 0$ and $\beta\in H_2(V)$ let $\Mb_{g,n}(V,\beta)$ denote
Kontsevich's moduli stack of $n$-pointed stable maps of genus $g$ and
class $\beta$ to $V$. 
For $g>1-n$ let
$\Mt_{g,n}$ denote one of the above stacks if $g=0$,
and $\Mt_{g,n}:=\Mb_{g,n}\times\Ws^n$ else.
In \S \ref{section stt}
we construct natural
morphisms 
$$
\stt:\Mb_{g,n}(V,\beta)\to \Mt_{g,n}
$$
such that in the stable range the composition 
of $\stt$ with the projection
$\pr_1:\Mt_{g,n}\to\Mb_{g,n}$ is the usual stabilisation
morphism $\sta:\Mb_{g,n}(V,\beta)\to\Mb_{g,n}$.
Recall that on $\Mb_{g,n}(V,\beta)$ there are $n$ natural
line bundles which can be defined as $x_i^*\omega$ where
$(\Cc\to V,(x_i))$ is the universal stable map over $\Mb_{g,n}(V,\beta)$,
and $\omega=\omega_{\Cc/\Mb}$
is the relative dualizing sheaf.
We show in Proposition \ref{stt*psi} that the Chern classes of
these line bundles are the pull-back of certain classes $\psit_i$
on $\Mt_{g,n}$.

Recall that by means of the morphism $\sta$ and the virtual fundamental class
on $\Mb_{g,n}(V,\beta)$ one can define the structure of an algebra
over the operad $(A_*\Mb_{0,n+1})_{n\geq 2}$ on the cohomology of $V$
with coefficients in some topological ring $\Lambda$.
In view of the analogy of $\stt$ to $\sta$ it seems reasonable
to conjecture 
that there is even a structure of an algebra
over the cyclic operad $(A_*\Mt_{0,n+1})_{n\geq 1}$ on $H^*(V,\Lambda)$.
Furthermore the existence of the classes $\psit_i$ makes it plausible
that gravitational correlators can be expressed in terms of that structure.
We confer to Conjectures \ref{conj1} and \ref{conj2} for more precise
statements.

As a first test for these conjectures we show in
\S\ref{section gravitation} 
how Theorem 1.2 of \cite{KM}
can be viewed as a consequence of 
a simple identity (cf. Lemma \ref{identity}) in $R[t_0]$.

As a second test we show how an
extension of the cited theorem (cf. Proposition \ref{formula})
to the unstable range can also be deduced from the conjectures.
The resulting formula can be proven by other means and
seems to be well-known.

Still assuming the conjectures we prove as a final result
(cf. Proposition \ref{final} and Remark \ref{finalrem})
that knowledge of the
$(A_*\Mt_{0,n+1})_{n\geq 1}$-algebra structure on $H^*(V,\Lambda)$
amounts to the same thing as the knowledge of
some numbers which probably can be expressed by means of gravitational
two-point correlators together gravitational
quantum cohomology in the stable range in the sense 
of \cite{Manin}, Ch VI, \S 2.2. 

This paper focuses attention on the geometry of the new cyclic operad
$\Mt_0$. It does not contain any new relation between gravitational
correlators. Nevertheless we hope that it provides a new point of
view on the geometric nature of these correlators, which might
be useful in the course of further research.

\vspace{2mm}
I would like to express my gratitude to Prof. Manin, who 
suggested to me that my insight in the geometry of modifications
of pointed curves might be useful for constructing an extended modular
operad.
Part of this work has been carried out during a visit at Max-Planck Institute
in Bonn and during a visit at the Alfred R\'enyi Mathematical Institute in Budapest.
The hospitality of these institutes is gratefully acknowledged.

\section{Modifications of pointed curves}
\label{section modifications}

Let $C$ be a curve and let $x\in C$ be a smooth point.
A {\em modification}  of $C$ at $x$ is a pointed curve $(C',x')$ together
with a morphism $f:C'\to C$ of pointed curves such that 
$f$ is an isomorphism outside $x$ and such that
the inverse image of $x$ is a chain of projective lines with
$x'$ in the extremal component. More precisely we require
that either $f$ is an isomorphism mapping $x'$ to $x$, or that
$C'\isomorph C\cup R_1\cup\dots\cup R_l$ for some $k\geq 1$,
where $R_i$ are 
projective lines such that 
\begin{itemize}
\item
$R_1\cap C=\{x\}$  and $x'\in R_l$,
\item
for $i=1,\dots,l-1$
the intersection
$R_i\cap R_{i+1}$ consists of one point $y_i$,
\item 
all the other intersections between the $C,R_1,\dots,R_l$
are empty,
\item
the points $x,y_1,\dots,y_{l-1},x'$ are distinct.
\end{itemize}
and $f$ maps all components $R_i$ to $x$.
If $f$ is an isomorphism, we say that $C'\to C$ is a modification
of length zero, otherwise we define the number $l$
to be the {\em length} of the modification.

The notion of modifications can readily be defined for families:
Let $S$ be a scheme, let $\pi:\Cc\to S$ be a flat family of curves
over $S$ and let $x:S\to\Cc$ be a section of $\pi$ whose image lies 
in the smooth locus of $\pi$.
A {\em modification} of $\Cc$ at $x$ is a flat family $\pi:\Cc'\to S$ of curves
together with a section $x':S\to \Cc'$ and an $S$-morphism $\Cc'\to \Cc$,
which geometric fibre-wise over $S$ is a modification in the sense above.
We say that the modification $\Cc'\to \Cc$ is of {\em length} $\leq l$,
if the length of the modifications at each geometric fibre is at most $l$.
Modifications of length $\leq 1$ will also be called {\em simple modifications}.

Let $f_1:(\Cc_1,x_1)\to(\Cc,x)$ and $f_2:(\Cc_2,x_2)\to(\Cc,x)$ be two modifications of
$\Cc$ at $x$. An {\em isomorphism} between these modifications is an isomorphism
$g:(\Cc_1,x_1)\to(\Cc_2,x_2)$ of pointed curves, such that $f_1=f_2\comp g$.
There is an obvious notion of pull-back for modifications:
Given a pointed curve $(\Cc,x)$ over a scheme $S$ and a modification $\Cc'\to \Cc$
of $\Cc$ at $x$, for any morphism $T\to S$ the fibre product 
$\Cc'\times_ST\to\Cc\times_ST$ is a modification 
of $\Cc_T:=\Cc\times_ST$ at the pull-back
section $x_T$.
It is easy to see that for a given pointed curve $(\Cc,x)$ over a scheme
$S$ the S-groupoid $\Tails(\Cc,x)$ 
which to each $S$-scheme $T$ associates the groupoid
of modifications of $\Cc_T$ at $x_T$ is a stack. 
We let $\Tails_l(\Cc,x)$ denote the substack of $\Tails(\Cc,x)$  which
parametrises modifications of length $\leq l$ of $\Cc$ at $x$.

By Proposition 5.3 in \cite{degeneration} we have:

\bpr
\label{T1 isom A1/Gm}
For any pointed curve
$(\Cc,x)$ over a scheme $S$ there is a canonical
isomorphism of groupoids
$$
\left\{
\begin{array}{lll}
\text{Pairs $(\Ll,\lambda)$ where $\Ll$}\\
\text{is a line bundle over $S$}\\
\text{and $\lambda$ is a section of $\Ll$}
\end{array}
\right\}
\to
\left\{
\begin{array}{ll}
\text{Simple modifications}\\
\text{of $\Cc$ at $x$}
\end{array}
\right\}
$$
which commutes with base change.
In other words, we have a canonical isomorphism
$$
[\A^1/\Gm]\times S\isomto\Tails_1(\Cc,x) 
\quad.
$$
\epr

As a matter of notation, for a pointed curve $(\Cc,x)$
over a scheme $S$, a line bundle $L$ over $S$ and
a global section $\lambda$ of $L$ we let
$$
(\Cc,x)\vdash(L,\lambda)
$$
denote the pointed curve $(\Cc',x')$, where $(\Cc',x')\to(\Cc,x)$
corresponds to $(L,\lambda)$ via the isomorphism in
Proposition \ref{T1 isom A1/Gm}.
We will often make use of the following result which is an immediate
consequence of the proof of Proposition 5.3 in \cite{degeneration}:

\blem
With the above notation let $(\Cc',x')=(\Cc,x)\vdash(L,\lambda)$.
Then we have
$$
(x')^*\Oo_{\Cc'}(-x')=x^*\Oo_{\Cc}(-x)\tensor L
\quad.
$$
\elem

Proposition \ref{T1 isom A1/Gm} implies
in particular that
the stack $\Tails_1(\Cc,x)$ is algebraic and has a very
concrete description which depends only on the base scheme,
not on the pointed curve $(\Cc,x)$.
In the next section we will prove an analogous result for $\Tails(\Cc,x)$
and for $\Tails_l(\Cc,x)$ for 
arbitrary $l$.
That result relies on the fact that locally over a base scheme
of finite type every modification of a proper curve can be decomposed 
into simple modifications. 
In the remaining of this section we prove this fact.

We begin with a definition.
If $f:(\Cc',x')\to(\Cc,x)$ is a modification of length $\leq l$ 
of $\Cc$ at $x$ and if
$g:(\Cc'',x'')\to(\Cc',x')$ is a modification of length $\leq l'$ of $\Cc'$ at $x'$ then
obviously the composition $f\comp g:(\Cc'',x'')\to(\Cc,x)$ is a modification 
of length $\leq l+l'$ of $\Cc$ at $x$.
A modification  $f:(\Cc',x')\to(\Cc,x)$ of $\Cc$ at $x$ is called {\em decomposable},
if it is the composition of modifications of simple modifications.

There are two important observations to make here.
The first one is that if a modification $f:(\Cc',x')\to(\Cc,x)$ is decomposable,
say 
$$
f:(\Cc',x')=(\Cc_l,x_l)\to(\Cc_{l-1},x_{l-1})\to\dots\to(\Cc_{1},x_{1})
\to(\Cc_{0},x_{0})=(\Cc,x)
$$ 
where each step is a simple modification, then 
the modifications $f_i:(\Cc_{i},x_{i})\to(\Cc_{i-1},x_{i-1})$ involved
are far from being uniquely determined 
and not even their number is determined by $f$.
This is clear since we are free to insert an isomorphism at any place.
But even if we require that none of the $f_i$ should be an isomorphism,
uniqueness is not guaranteed. 
Indeed let $S=\A^1$, $\Cc:=\Pp^1_S$ and $x:S\to\Cc$ any (constant, say)
section. Let $\Cc'$ be the blow up of $\Cc$ at the closed subset 
consisting of the union of the points $x(0)$ and $x(1)$, 
and let $x'$ be the proper transform of $x$. 
Then $f:\Cc'\to\Cc$ itself is simple, but there exist also two
different representations of $f$ as the composition of two simple
modifications,
since $\Cc'$ may be obtained by either first
blowing up $x(0)$ and then $x(1)$, or vice versa.

The second important observation is that decompositions don't necessarily exist.
Indeed,
let $S$ be an algebraic surface (e.g. $S=\Pp^2$), 
let $D_0,D_1,D_2\subset S$ be curves
on $S$ with pairwise transversal intersection and empty triple
intersection. Let $U_i:=S\setminus D_i$ ($i\in\Z/3$) and
let $X_i$ be the blowing up of $\Pp^1\times U_i$ first along
$\{0\}\times(D_{i+1}\cap U_i)$ and then along the proper
transform of $\{0\}\times(D_{i+2}\cap U_i)$.
Then there are canonical isomorphisms 
$
X_i\times_{U_i}(U_i\cap U_j)\isomorph X_j\times_{U_j}(U_i\cap U_j)
$
by which we can glue together the $X_i$ to obtain a curve $X$
over $S$ whose schematic picture reminds one of M. C. Escher's illusion
of an infinitely ascending staircase:
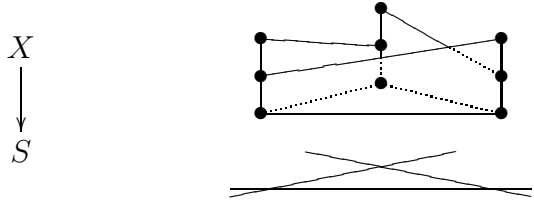
\[
\label{escher}
\vcenter{
\xymatrix{
X \ar[d] \\
S
}
}
\qquad\qquad\qquad
\vcenter{
\xy
0;<1cm,0cm>:
(-2,0); (2,0)
**@{-};
(-2,-0.1); (1,0.5)
**@{-};
(2,-0.1); (-1,0.5)
**@{-};
(-1.6,1); (1.6,1)
**@{-};
(-1.6,1); (0,1.4)
**@{..};
(1.6,1); (0,1.4)
**@{..};
(-1.6,1); (-1.6,2)
**@{-};
(1.6,1); (1.6,2)
**@{-};
(0,1.4); (0,1.75)
**@{..};
(0,1.75); (0,2.4)
**@{-};
;
(-1.6,1)*{\bullet};
(1.6,1)*{\bullet};
(0,1.4)*{\bullet};
;
(-1.6,1.5)*{\bullet};
(1.6,1.5)*{\bullet};
(0,1.9)*{\bullet};
;
(-1.6,2)*{\bullet};
(1.6,2)*{\bullet};
(0,2.4)*{\bullet};
;
(-1.6,1.5); (1.6,2) **@{-};
(1.6,1.5); (0.912,1.887) **@{..}; 
(0.912,1.887); (0,2.4) **@{-};
(0,1.9); (-1.6,2) **@{-};
\endxy
}
\]
Let $x:S\to X$ be the section of $X\to S$ which over $U_i$
is the proper transform of the zero section of $\Pp^1_{U_i}$.
The construction of $X$ provides us with a map $f:(X,x)\to(\Pp^1_S,0)$
which clearly is a modification of length $\leq 2$ of $\Pp^1_S$ at $0$.
Suppose there was a decomposition $f=f_1\comp f_2$ of $f$ into 
simple modifications. Then over $U_0$ the decomposition
is necessarily the successive blow down given by the construction.
This means in particular that $f_1$ is an isomorphism over 
$U_0\cap U_2$ but not over $U_0\cap U_1$. Applying the same
argument to $U_1$ we see that $f_1$ is an isomorphism over
$U_1\cap U_0$ but not over $U_1\cap U_2$. These requirements
clearly cannot be satisfied simultaneously.

The following Proposition gives a criterion for the decomposability
of modifications. At the same time it shows that the morphism $X\to S$
(Escher's staircase) is not projective.
\begin{proposition}
Let $(\Cc,x)$ be a pointed curve over a 
quasi-compact scheme $S$  and let 
$f:(\Cc',x')\to(\Cc,x)$ be a modification of $\Cc$ at $x$.
Suppose there exists a relatively (over $S$)
ample line bundle $E$ on $\Cc'$. Then $\Cc'\to\Cc$ is decomposable.
\end{proposition}

\bpf
Replacing $E$ by a suitable power, we may assume that $E$ is 
relatively very ample.
We define maps $n_{\Cc'}, l_{\Cc'}:S\to\N_0$ as follows.
If $\Cc'_z\to\Cc_z$ is an isomorphism we let $n_{\Cc'}(z)=l_{\Cc'}(z)=0$. 
Otherwise let $n(z)\in \N$ be the degree of
$E$ restricted to the component $R(z)$ of the fibre $\Cc'_z$ which
contains $x(z)$, and let $l_{\Cc'}(z)$ be the length of the modification 
$\Cc'_z\to\Cc_z$.
Then $n_{\Cc'}$ and $l_{\Cc'}$ are constructible
functions. Let $N$ be the maximum of $n_{\Cc'}$.
The line bundle $E\tensor_{\Oo}\Oo(-Nx)$ is generated 
by global sections and thus defines a morphism $f_1$ into some
high dimensional projective space which contracts precisely
those $R(z)$ for which $n(z)$ is maximal. 
Let $\Cc_1$ be the image of $f_1$ and let $x_1:=f_1\comp x$.
It is clear that $f_1:(\Cc',x')\to(\Cc_1,x_1)$
is a simple modification and that $f$ factorises over $f_1$.
Furthermore, by construction there is a relatively ample line bundle $E_1$
on $\Cc_1$ and the function $l_{\Cc_1}$ is strictly smaller than $l_{\Cc'}$.
Therefore iterating this procedure we obtain the sought for decomposition.
\epf

\bco
\label{loc dec}
Let $(\Cc,x)$ be a proper pointed curve over a 
scheme $S$  and let 
$f:(\Cc',x')\to(\Cc,x)$ be a modification of $\Cc$ at $x$.
Then locally over $S$ the modification is decomposable.
\eco

\bpf
This is immediate from the above Proposition, since every proper
curve is locally projective over the base.
\epf

\section{The stack $\Ws$}
\label{section W}

In this section we prove that for a given pointed curve 
$(\Cc,x)$ over some base scheme $B$ the stack of modifications $\Tails(\Cc,x)$ 
introduced in the previous paragraph is naturally isomorphic
to $\Ws\times B$ for a certain stack $\Ws$ which has a very
concrete description.
 
To define $\Ws$ we proceed in three steps. First
we define a $\C$-groupoid $\Ws''$, which is not even a prestack.
By sheafifying the presheaf of isomorphisms in $\Ws''$ we will obtain
a prestack $\Ws'$ which is not a stack. Finally, by forcing descent
data to be effective, we obtain the stack $\Ws$. 
The passage from $\Ws''$ to $\Ws$ follows the general prescription
as given in \cite{L-MB} (3.2.1)(3).

In fact, for each of the groupoids $\Ws''$, $\Ws'$, $\Ws$ and each non-negative
integer $l$ there is a level-$l$-version $\Ws''_l$, $\Ws'_l$, $\Ws_l$ and
we have
$$
\Ws=\underset{l}{\varinjlim}\Ws_l
$$
and likewise for $\Ws''$ and $\Ws'$.

We will show that
there is a representable smooth surjective morphism
$$
[\A^1/\Gm]^l \to \Ws_l
$$
which implies in particular that $\Ws_l$ is algebraic.
This will prove that $\Ws$ is an ind-object in the category 
of algebraic stacks.

\vspace{5mm}
{\bf Definition of $\Ws''$:}
Let $U\in\ob(\Aff/\C)$. We define the groupoid $\Ws''(U)$ as follows.
The objects of $\Ws''(U)$ are finite (possibly empty) sequences
$$
(L_1,\lambda_1)(L_2,\lambda_2)\dots(L_r,\lambda_r)
$$
where $L_i$ is an invertible $\Oo_U$-module and $\lambda_i$ is
a section of $L_i$ which is {\em not invertible} in the sense
that $\lambda_i$ vanishes in at least one point of $U$.
The empty sequence will be denoted by $\emptyset$.
A morphism
$$
\alpha:
(L_1,\lambda_1)\dots(L_r,\lambda_r)
\to
(L'_1,\lambda'_1)\dots(L'_r,\lambda'_{r'})
$$
in $\Ws''(U)$ exists only if $r=r'$.
If $r=r'\geq 1$, then $\varphi$ is simply
a sequence of isomorphisms $\alpha_i:L_i\isomto L'_i$ such
that $\alpha_i(\lambda_i)=\lambda'_i$.
If $r=r'=0$ then both objects are empty and there is only
the identity morphism. 

To each arrow $\varphi:V\to U$ in $(\Aff/\C)$ we associate
the functor $\varphi^*:\Ws''(U)\to\Ws''(V)$ which
to an object
$$
(L_1,\lambda_1)\dots(L_r,\lambda_r)
$$
in $\Ws''(U)$
associates the object
$$
(\varphi^*L_{i_1},\varphi^*\lambda_{i_1})
\dots
(\varphi^*L_{i_s},\varphi^*\lambda_{i_s})
$$
in $\Ws''(V)$, where $1\leq i_s<\dots i_s\leq r$ is the (possibly empty)
set of all indices $j\in[1,r]$ such that $\varphi^*\lambda_j$ is not
invertible. The pull back of a morphism in $\Ws''(U)$ is the obvious
one.

Clearly, if $U$ is an object in $(\Aff/\C)$ and $\xi$, $\eta$ are
objects in $\Ws''(U)$, then the presheaf
$$
\Isom_{\Ws''}(\xi,\eta):
\left\{
\begin{array}{ll}
(\Aff/U) & \to (\Ens) \\
(V\to U) & \mapsto \Hom_{\Ws''(V)}(\xi_V,\eta_V)
\end{array}
\right.
$$
will in general {\em not} be a sheaf.
To illustrate this, let $\xi=(L_1,\lambda_1)(L_2,\lambda_2)$
and let $\eta=(L_2,\lambda_2)(L_1,\lambda_1)$, where
the closed subsets $Z_1:=\{\lambda_1=0\}$ and $Z_2:=\{\lambda_2=0\}$ of
$U$ are non-empty and disjoint.
Then if we put $U_i:=U\setminus Z_i$, we have
\begin{eqnarray*}
\xi_{U_1}= & (L_2,\lambda_2)|_{U_1} & =\eta_{U_1}
\quad\text{and}
\\
\xi_{U_2}= & (L_1,\lambda_1)|_{U_2} & =\eta_{U_2}
\end{eqnarray*}
but the set $\Hom_{\Ws''(U)}(\xi,\eta)$ is empty.
Thus $\Ws''$ is not a prestack in the sense of
\cite{L-MB} (3.1).

\vspace{5mm}
{\bf Definition of $\Ws'$:}
Let $U$ be an object in $(\Aff/\C)$. 
We define $\Ws'(U)$ as the category which has the
same objects as $\Ws''(U)$ but whose set of 
morphisms
$\Hom_{\Ws'(U)}(\xi,\eta)$
for $\xi,\eta\in\ob(\Ws'(U))=\ob(\Ws''(U))$ is the set of
global sections of the sheafification of the
(Zariski-) presheaf 
\[
\Isom_{\Ws''}^{\Zar}(\xi,\eta):
\left\{
\begin{array}{ll}
\{\text{Zariski open sets in $U$}\} & \to (\Ens) \\
(V\subseteq U) & \mapsto \Hom_{\Ws''(V)}(\xi_V,\eta_V)
\end{array}
\right.
\]
In more concrete terms, let
\begin{eqnarray*}
\xi &=&
(L_1,\lambda_1)\dots(L_r,\lambda_r)
\\
\eta &=&
(M_1,\mu_1)\dots(M_s,\mu_s)
\end{eqnarray*}
and for each $x\in U$ let 
$I(x)\subseteq[1,r]$ 
($J(x)\subseteq[1,s]$)
be the ordered set of indices $i$ with $\lambda_i(x)=0$
(with $\mu_i(x)=0$).
Then we have 
\begin{eqnarray*}
\xi_x &=& ((L_{i,x},\lambda_{i,x}), i\in I(x))
\\
\eta_x &=& ((M_{i,x},\mu_{i,x}), i\in J(x))
\end{eqnarray*}
An arrow $\xi\to\eta$ in the category $\Ws'(U)$
is a mapping 
\[
\alpha:
U\to\coprod_{x\in U}\Hom(\xi_x,\eta_x)
\]
such that $\alpha(x)\in\Hom(\xi_x,\eta_x)$ and 
for every $x\in U$ there is a Zariski open neighbourhood $V\subseteq U$
of $x$ and a morphism $\beta\in\Hom_{\Ws''(V)}(\xi_V,\eta_V)$
such that $\alpha(y)=\beta_y$ for all $y\in V$.

For an arrow $\varphi:V\to U$ in $(\Aff/\C)$ the functor
$\varphi^*:\Ws'(U)\to\Ws'(V)$ is the same as in the case of $\Ws''$ 
at least on the level of objects. On the level of morphisms this
functor is the obvious one.

We will now give an example which shows that $\Ws'$ is not a stack.
Let $U$ be an algebraic surface (one could take $U=\A^2$) and let
$Z_0$, $Z_1$, $Z_2$ be three curves on $U$ which intersect pairwise
in one point but whose triple intersection is empty.
Let $U_i\subset U$ be the complement of $Z_i$ in $U$.
For $i\in\Z/3$ let $L_i$ be an invertible sheaf on $U$
which has a section $\lambda_i$ whose set of zeroes is
precisely $Z_i$.
For $i\in\Z/3$  we define the following object
\[
\xi_i:=
(L_{i+1}|_{U_i},\lambda_{i+1}|_{U_i})
(L_{i+2}|_{U_i},\lambda_{i+2}|_{U_i})
\]
in $\Ws'(U_i)$.
Then we have natural isomorphisms
\[
f_i:
\xi_i|_{U_i\cap U_{i+1}}\isomto
(L_{i+2}|_{U_i\cap U_{i+1}},\lambda_{i+2}|_{U_i\cap U_{i+1}})\isomto
\xi_{i+1}|_{U_i\cap U_{i+1}}
\]
in $\Ws'(U_i\cap U_{i+1})$.
Since for all $i\in\Z/3$ the object
$
\xi_i|_{U_0\cap U_1\cap U_2}
$ 
is the empty one
in $\Ws'(U_0\cap U_1\cap U_2)$,
the data $(\xi_i,f_{i,j})$ form a descent datum relative to the
covering $(U_i\injto U)_i$.
It is easy to see that this descent datum is not effective.

The reader will have noticed the similarity of this example
with ``Escher's staircase'' from the previous paragraph.
In fact the two examples correspond via the isomorphism
$\Tails(\Pp^1,x)\isomto\Ws\times U$ of Proposition \ref{BxW isom T}
below.

\vspace{5mm}
{\bf Definition of $\Ws$:}
Let $\Ws$ be the stack associated to the 
prestack $\Ws'$.
(cf. \cite{L-MB}, lemme (3.2)).

For $U\in\ob(\Aff/\C)$ the category $\Ws(U)$ has as
objects triples 
\[
(U'\to U,\xi',f'')
\quad,
\]
where $U'\to U$ is a covering with one element 
in $(\Aff/\C)$, $\xi'$ is an
object in $\Ws'(U')$
and $f''$ is a descent datum for $\xi'$
with respect to the covering $U'\to U$.

\blem
\label{zar}
One obtains an isomorphic stack, if one takes only triples 
$(U'\to U,\xi,f'')$, where
$U'\to U$ is a {\em Zariski} covering.
\elem

\vspace{5mm}
{\bf Definition of 
$\Ws''_l$,
$\Ws'_l$, and
$\Ws_l$
:}
For $l\geq 0$ we let $\Ws''_l$ be the subgroupoid of $\Ws''$ consisting
of sequences
\[
(L_1,\lambda_1)(L_2,\lambda_2)\dots(L_r,\lambda_r)
\]
with $r\leq l$.
The prestack $\Ws'_l$ (the stack $\Ws_l$) is constructed 
from $\Ws''_l$ (from $\Ws'_l$) analogously
as  $\Ws'$ from $\Ws''$ (as $\Ws$ from $\Ws'$).
We have obvious inclusion morphisms
\[
\Ws''_0\injto\Ws''_1\injto\dots\injto\Ws''_{l-1}\injto\Ws''_l\injto\dots
\injto\Ws''
\]
and analogously for $\Ws'$ and $\Ws$.
It is clear that
\[
\Ws''=\underset{l}{\varinjlim}\Ws''_l
\quad, \quad
\Ws'=\underset{l}{\varinjlim}\Ws'_l
\quad, \quad
\Ws=\underset{l}{\varinjlim}\Ws_l
\quad.
\]
Furthermore, we have natural morphisms
\[
\Ws''_l\to\Ws'_l\to\Ws_l
\quad\text{and}\quad
\Ws''\to\Ws'\to\Ws
\]
compatible with the inclusion morphisms.

\begin{proposition}
\label{etale}
Let $l$ be a non-negative integer and let
\[
[\A^1/\Gm]^l \to \Ws''_l
\]
be the morphism of $\C$-groupoids which sends an object 
\[
((L_1,\lambda_1),\dots,(L_l,\lambda_l))
\]
of $[\A^1/\Gm]^l(U)$ to the object
\[
(L_{i_1},\lambda_{i_1})\dots(L_{i_r},\lambda_{i_r})
\] 
of $\Ws''_l(U)$, where $1\leq i_1<\dots<i_r\leq l$
is the set of all indices $j\in[1,l]$ such that the section $\lambda_j$ 
is not invertible.
Then the composite morphism
\[
[\A^1/\Gm]^l \to \Ws''_l \to \Ws'_l \to \Ws_l
\]
is a representable surjective \'etale morphism between stacks.
\end{proposition}

\begin{proof}
Let $U$ be an affine $\C$-scheme
and let $\xi\in\Ws_l(U)$. By \ref{zar}
for each $x\in U$ there is an open neighbourhood 
$U_x$ of $x$ such that
$$
\xi|_{U_x}=(L_{x,1},\lambda_{x,1})\dots(L_{x,m_x},\lambda_{x,m_x})
$$
for some line bundles $L_{x,i}$ over $U_x$ with 
sections $\lambda_{x,i}$ such that $m_x\leq l$ and
$\lambda_{x,i}(x)=0$ for all $i\in[1,m_x]$.
Let 
$$
X':=\coprod_{x\in U}\ \coprod_{\alpha:[1,m_x]\injto[1,l]} U_{x,\alpha}
$$
where 
$U_{x,\alpha}$ is a copy of $U_x$, and 
$\alpha$ runs through all strictly increasing maps 
$[1,m_x]\to[1,l]$.
We define an equivalence relation $\sim$ on $X'$ by
requiring $u\sim v$ for $u\in U_{x,\alpha}$, $v\in U_{y,\beta}$
if and only if $u$ and $v$ coincide as points of U, and
if furthermore the equality
$$
\alpha\left(\{i\in[1,m_x]\st \lambda_{x,i}(u)=0\}\right)=
\beta\left(\{i\in[1,m_y]\st \lambda_{y,i}(v)=0\}\right)
$$
of subsets of $[1,l]$ holds. Since the set of all points of
$U_{x,\alpha}$ which are equivalent to some point in $U_{y,\beta}$
is open, the quotient space
$$
X:=X'/\sim
$$
has a natural (in general non-separated) scheme structure. 
Furthermore, the open immersions
$U_{x,\alpha}\to U$ induce a surjective morphism $\pi:X\to U$ which
is a Zariski-local isomorphism and hence \'etale.

Over $U_{x,\alpha}$ we have an $l$-tuple 
$
((M_{x,1},\mu_{x,1}),\dots,(M_{x,l},\mu_{x,l}))
$
of line bundles with sections where
$$
(M_{x,j},\mu_{x,j})=
\bca
(L_{x,i},\lambda_{x,i})
& \text{if there is an $i\in[1,m_x]$ with $\alpha(i)=j$,}\\
(\Oo_{U_{x,\alpha}},1) & \text{else.}
\eca
$$
and the resulting $l$-tuple of line bundles with sections
on $X'$ descends to an $l$-tuple 
$$
\eta:=((M_1,\mu_1),\dots,(M_l,\mu_l))
$$
on $X$, which can be viewed as a morphism $\eta:X\to[\A^1/\Gm]^l$.
We claim that the diagram
$$
\xymatrix{
X \ar[r]^(.4)\eta \ar[d]_\pi 
&
\text{$[\A^1/\Gm]^l$} \ar[d]^h
\\
U \ar[r]^\xi
&
\text{$\Ws_l$}
}
$$
is Cartesian.

Indeed, let $V$ be a scheme, $f:V\to U$ a morphism and
$\tau=((N_1,\nu_1),\dots,(N_l,\nu_l))$ an $l$-tuple of line
bundles with sections on $V$ and assume that we are given an 
isomorphism $\varphi:f^*\xi\isomto h(\tau)$. 
We have to show that there is a unique
morphism $g:V\to X$ such that (i) $f=\pi\comp g$,
(ii) there is an isomorphism $\psi:\tau\isomto g^*\eta$, and 
(iii) the pull back via $g$ of the isomorphism $\pi^*\xi\isomto h(\eta)$
is $h(\psi)\comp\varphi$.
Let $p$ be a point of $V$ and let $x=f(p)$. 
It is clearly sufficient to
prove the existence of a unique morphism $V'\to X$ with
the above properties, where $V'\subset V$ is some open neighbourhood
of $p$. Let $I\subset[1,l]$ be the set
of all indices $i$ with $\nu_i(p)=0$ and let $\alpha:[1,n]\injto[1,l]$
be the unique strictly increasing map with image $I$.
Observe that we have $n=m_x$ since $f^*\xi\isomorph h(\tau)$.
Replacing $V$ by an open neighbourhood of $p$, we may assume that
$f$ factorises over $U_x$ and that $\varphi$ is given by the
$n$ isomorphisms 
$\varphi_i:
f^*(L_{x,i},\lambda_{x,i})\isomto(N_{\alpha(i)},\nu_{\alpha(i)})$ 
($i=1,\dots,n$). We define $g$ as the composition
$$
V\overset{f}{\to}U_x=U_{x,\alpha}\to X
\quad.
$$
Property (i) for $g$ then holds trivially.
To check property (ii) for $g$ we first observe that
by definition we have $g^*\eta=((M'_1,\mu'_1),\dots,(M'_l,\mu'_l))$
where $(M'_j,\mu'_j)=f^*(L_{x,i},\lambda_{x,i})$ if there
is an $i$ with $j=\alpha(i)$, and $(M'_j,\mu'_j)=(\Oo_V,1)$ else.
We define $\psi:=(\psi_1,\dots,\psi_l)$ where $\psi_j:=\varphi_i^{-1}$
if there is an $i$ with $j=\alpha(i)$.
If there is no such $i$ then $\psi_j$ is the
unique isomorphism $N_j\to\Oo_V$ of line bundles, which sends the
invertible section $\nu_j$ to $1$. Property (iii) for $g$
is now immediate and it is also clear that $g$ is uniquely determined
by  properties (i)-(iii).

This proves the first part of the proposition.
It implies in particular that the composed morphism
$$
\A^l\to[\A^1/\Gm]^l\to\Ws_l
$$
is representable, surjective, and smooth.
Therefore for the second part of the proposition
all that remains to be shown is
that for every affine $\C$-scheme $U$ and objects
$\xi$, $\eta$ of $\Ws_l(U)$ the sheaf $\Isoms(\xi,\eta)$
of isomorphisms of $\xi$ to $\eta$ is representable
by a scheme $\Isom(\xi,\eta)$ which is separated 
quasi-compact over $U$.

For this we first consider the case where $\xi=(L,\lambda)$
and $\eta=(M,\mu)$ for some line bundles with sections
$(L,\lambda)$ and $(M,\mu)$ on $U$.
In this case $\Isoms(\xi,\eta)$ is representable by a
locally closed subscheme of $X:=\Spec(\Sym^{\bullet}(L\tensor M^\dual))$.
More precisely we have $\Isom(\xi,\eta)=X^o\cap Y$ where
$X^o\subset X$ is the complement of the zero section of $X\to U$
and $Y\subset X$ is the closed subscheme corresponding to the
ideal generated by
$$
\{\omega(\mu)-\lambda\tensor\omega\st \omega\in M^\dual\}
\ \subset\ \Sym^\bullet(L\tensor M^\dual)
\quad.
$$
If we further specialise to the case where $L\isomorph\Oo_U\isomorph M$
then we find 
$$
\Isom(\xi,\eta)\isomorph\Spec(\Oo_U[T,T^{-1}]/(\mu-\lambda T))
$$
which shows in particular that $U'\times_U\Isom(\xi,\eta)=U'$ 
for open $U'\subset U$ over which both $\lambda$ and $\mu$ 
are invertible.

As a next step we assume that
$$
\xi=(L_1,\lambda_1)\dots(L_n,\lambda_n)
\qquad,\qquad
\eta=(M_1,\lambda_1)\dots(M_n,\lambda_n)
$$
and that $(L_i,\lambda_i)\isomorph(M_i,\mu_i)$ for  
$i=1,\dots,n$.
In this case we have
$$
\Isom(\xi,\eta)=\prod_{i=1}^n\Isom((L_i,\lambda_i),(M_i,\mu_i))
\quad.
$$

For general $\xi$, $\eta$ let $U'$ be the open subset of $U$
which consists of all points $x$ for which the germs 
$\xi_x$ and $\eta_x$ are isomorphic, and let
$U'=\cup_{i\in I}U_i$ be an open affine covering 
such that 
$$
\xi|_{U_i}=(L_{i,1},\lambda_{i,1})\dots(L_{i,n_i},\lambda_{i,n_i})
\qquad,\qquad
\eta|_{U_i}=(M_{i,1},\lambda_{i,1})\dots(M_{i,n_i},\lambda_{i,n_i})
$$
for some line bundles with sections 
$(L_{i,j},\lambda_{i,j})\isomorph(M_{i,j},\mu_{i,j})$ over $U_i$.
Let $X_i:=\Isom(\xi|_{U_i},\eta|_{U_i})$. 
Then we have natural isomorphisms
$$
(U_i\cap U_j)\times_{U_i}X_i\isomorph
\Isom(\xi|_{U_i\cap U_j},\eta|_{U_i\cap U_j})\isomorph
(U_i\cap U_j)\times_{U_j}X_j
$$
by which we can glue the $X_i$ together to form the scheme
$\Isom(\xi,\eta)$. It is clear from the construction that
it is separated and quasi-compact over $U$.
\end{proof}

As an immediate consequence of Proposition \ref{etale} we have:
\bco
\label{etale cor}
The groupoids $\Ws$ and $\Ws_l$ are smooth algebraic stacks
of dimension $0$. Moreover the $\Ws_l$ are of finite type
and we have the following inclusion of open substacks
$$
\Ws_1\subset\Ws_2\subset\dots\subset\Ws_l\subset\dots\subset\Ws
\quad.
$$
\eco

The following result explains the relationship between the stack
$\Ws$ and the groupoid of modifications defined in the previous paragraph.
\bpr
\label{BxW isom T}
Let $B$ be a $\C$-scheme and let $(\Cc,x)$ be a pointed curve
over $B$.
Then there is a natural isomorphism of $B$-stacks
\[
B\times\Ws\isomto\Tails(\Cc,x)
\quad.
\]
For each $l\geq 0$ it maps the open substack $\Ws_l$
isomorphically onto $\Tails_l(\Cc,x)$.
\epr

\bpf
Let $S$ be a $B$-scheme and
let $\xi=(L_1,\lambda_1)\dots(L_r,\lambda_r)$ be an object in
$\Ws''(S)$. Repeatedly applying Proposition \ref{T1 isom A1/Gm} we obtain
a sequence
$$
(\Cc',x'):=(\Cc_r,x_r)\to\dots\to(\Cc_1,x_1)\to(\Cc_0,x_0):=(\Cc_S,x_S)
$$
of simple modifications. In this sequence the modification
$(\Cc_i,x_i)\to(\Cc_{i-1},x_{i-1})$
corresponds to the pair $(L_i,\lambda_i)$.
The association $\xi\mapsto(\Cc',x')$ defines a morphism 
$\Ws''(S)\to\Tails(\Cc,x)(S)$ of groupoids.
Observe that if $T\to S$ is a morphism of $B$-schemes then the diagram
$$
\xymatrix{
\text{$\Ws''(S)$} \ar[r]\ar[d] & \text{$\Tails(\Cc,x)(S)$} \ar[d]\\
\text{$\Ws''(T)$} \ar[r] & \text{$\Tails(\Cc,x)(T)$}
}
$$
is 2-commutative. 
Indeed, this is a consequence of the fact that the isomorphism
of Proposition \ref{T1 isom A1/Gm} sends a pair $(L,\lambda)$ to
a modification of length zero if $\lambda$ is an invertible section of $L$.
Thus we see that there is a natural morphism
$
\Ws''\times B\to\Tails(\Cc,x)
$
of $B$-groupoids. Since $\Tails(\Cc,x)$ is a stack, this induces a morphism
$$
\Ws\times B\to\Tails(\Cc,x)
\quad.
$$

To construct an arrow in the reverse direction, we first 
consider the following situation.
Let $S$ be a $B$-scheme and let $(\Cc',x')\to(\Cc_S,x_S)$
be a modification which has a decomposition 
$$
(\Cc',x')=(\Cc_r,x_r)\to\dots\to(\Cc_1,x_1)\to(\Cc_0,x_0):=(\Cc_S,x_S)
$$
into simple modifications. By Proposition \ref{T1 isom A1/Gm}
to each step $(\Cc_i,x_i)\to(\Cc_{i-1}\to(\Cc_{i-1},x_{i-1})$
there is associated a pair $(L_i,\lambda_i)$ consisting of a
line bundle and a section. Thus we obtain an object
$$
\xi=(L_{i_1},\lambda_{i_1})\dots(L_{i_k},\lambda_{i_k})
\in \Ws''(S)
$$
where $i_1<\dots<i_k$ is the set of all indices $j\in[1,r]$ 
such that $\lambda_j$ is not invertible.
Furthermore, if $T\to S$ is a morphism of $B$-schemes and we apply
the above construction to the decomposition of $(\Cc'_T,x'_T)\to(\Cc_T,x_T)$
obtained by pulling back the given decomposition over $S$, then 
the object $\eta\in\Ws''(T)$ obtained in this way is canonically
isomorphic to the restriction $\xi_T$.

As we have seen at the end of \S \ref{section modifications},
a decomposition of a modification is not uniquely determined.
Thus also the element $\xi$ constructed above is ambiguous.
In fact the example in \S \ref{section modifications} shows that
two different decompositions may lead to non-isomorphic objects
in $\Ws''(S)$. We claim however that $\xi$ does not depend
on the given decomposition up to a canonical isomorphism
in $\Ws'(S)$.

To prove the claim we observe that in the case where $S$ is
the spectrum of a local ring, the decomposition of a modification
$(\Cc',x')\to(\Cc_S,x_S)$ is determined up to a unique isomorphism,
if one requires that none of the steps in the decomposition
is a modification of length zero. Indeed in this case each step corresponds
to the contraction of a rational component in the central fibre.
From this observation it follows that 
two decompositions of a given modification have isomorphic
germs at each point, after omitting the steps of length zero,
and that furthermore the isomorphisms between the germs are unique.
Thus between the germs of the objects in $\Ws''(S)$, which correspond
to the two decompositions, there are canonical isomorphisms which
fit together to form an isomorphism in $\Ws'(S)$.

By Corollary \ref{loc dec}, for any modification
$(\Cc',x')\to(\Cc_S,x_S)$ over a $B$-scheme $S$
there is an open covering $S=\cup S_i$ such that the modification 
$(\Cc'_{S_i},x'_{S_i})\to(\Cc_{S_i},x_{S_i})$ is decomposable for every $i$.
By the construction above we can associate an object
$\xi_i\in\Ws'(S_i)$ to this modification.
Let $S_{i,j}:=S_i\cap S_j$. 
Since both $\xi_i|_{S_{i,j}}$ and $\xi_j|_{S_{i,j}}$
are isomorphic to the object in $\Ws'(S'')$ which
we obtain by applying the above construction to the modification 
$(\Cc'_{S_{i,j}},x'_{S_{i,j}})\to(\Cc_{S_{i,j}},x_{S_{i,j}})$ we see that
there is a canonical isomorphism
$f_{i,j}:\xi_i|_{S_{i,j}}\isomto\xi_j|_{S_{i,j}}$.
Thus there is a uniquely determined 
object $\xi\in\Ws(S)$ with $\xi|_{S_i}=\xi_i$ for every $i$.

The association $((\Cc',x')\to(\Cc_S,x_S))\mapsto\xi$
defines a morphism 
$$
\Tails(\Cc,x)\to\Ws\times B
\quad.
$$
which is inverse to the morphism constructed at the beginning
of the proof.

The last part of the proposition is clear by construction.
\epf

\section{Intersection theory for exhaustible stacks}
\label{section exhaustible}

For an algebraic stack $\Xs$ of finite type over
$\C$, let $A_k\Xs$ be the $k$-th  homology as defined by 
Kresch (cf. \cite{Kresch}, Def. 2.1.11), but with coefficients in $\Q$.
In the following chapters we need an extension of the homology
functor to stacks like $\Ws$, which are not of finite type
but are {\em exhaustible} in the following sense.

\bdefin
Let $\Xs$ be an algebraic stack.
\balph
\item
An {\em exhaustion} for $\Xs$ is increasing series of open substacks of
finite type
$$
\Xs_1\subset\Xs_2\subset\dots\subset\Xs_n\subset\dots\subset\Xs
$$
such that $\Xs=\cup_n\Xs_n$.
\item
The stack $\Xs$ is called {\em exhaustible} if there is an
exhaustion for $\Xs$.
\item
Two exhaustions $(\Xs_n\subset\Xs)_n$ and $(\Xs'_n\subset\Xs)_n$
are called {\em compatible} if for each $n$ there is
an $m$ with $\Xs_n\subset \Xs'_m$, and for each $k$ there is 
an $l$ with $\Xs'_k\subset\Xs_l$.
\item
Let $f:\Xs\to\Ys$ be a morphism between two exhaustible algebraic
stacks and let  $(\Xs_n\subset\Xs)_n$ and $(\Ys_n\subset\Ys)_n$
be exhaustions for $\Xs$ and $\Ys$ respectively. We call $f$
to be {\em compatible} with the given exhaustions if 
$(f^\mo(\Ys_n)\subset\Xs)_n$ is an exhaustion for $\Xs$ which
is compatible with $(\Xs_n\subset\Xs)_n$.
\ealph
\edefin

We gather some facts about exhaustible stacks and their homology
theory. The proofs are straightforward and left to the reader.

\blem
Let $f:\Xs_1\to\Xs_2$ and $g:\Xs_2\to\Xs_3$ be two morphisms
between exhaustible algebraic stacks, and assume $f,g$ to be
compatible with the exhaustions. 
Then the composition $g\comp f$ is also compatible with the
exhaustions.
\elem

\blem
\label{exhaust Cart}
Assume that we have a Cartesian diagram
$$
\xymatrix{
\Xs'\ar[r]^{f'}\ar[d]_{g'}
&
\Ys'\ar[d]^g
\\
\Xs\ar[r]^f
&
\Ys
}
$$
where $\Xs,\Ys,\Ys'$ are exhaustible with respective exhaustions
$(\Xs_n\subset\Xs)_n$, $(\Ys_n\subset\Ys)_n$, 
$(\Ys'_n\subset\Ys')_n$.
Assume furthermore that $f$ and $g$ are compatible with the given
exhaustions. Then $\Xs'$ is exhaustible with exhaustion
$(f^\mo(\Ys_n)\times_{\Ys_n}g^\mo(\Ys_n)\subset\Xs')_n$,
and the morphisms $f'$ and $g'$ are compatible with respect
to the given exhaustions.
\elem

\bdefin
\label{A_* exhaust}
Let $\Xs$ be an exhaustible algebraic stack 
and let $(\Xs_n\subset\Xs)_n$ be an exhaustion of $\Xs$.
We define $A_i((\Xs_n\subset\Xs)_n):=\varprojlim_nA_i(\Xs_n)$
where the limit is taken with respect to the pull back morphisms
associated to the inclusions $\Xs_n\injto\Xs_{n+1}$.
Let $A_*((\Xs_n\subset\Xs)_n):=\Oplus_iA_i(\Xs)$.
\edefin

\brem
We have defined $A_*((\Xs_n\subset\Xs)_n)$ as the direct sum of the projective
limits $\varprojlim_nA_i(\Xs_n)$. We could equivalently (and ore concisely) 
have defined $A_*((\Xs_n\subset\Xs)_n):=\varprojlim'_nA_*(\Xs_n)$ where $\varprojlim'$
is the projective limit in the category of graded groups.
\erem

\blem
\label{lem exhaustion}
Let $\Xs$ be an exhaustible
algebraic stack and let $(\Xs_n\subset\Xs)_n$
be an exhaustion.
Then the group $A_i((\Xs_n\subset\Xs)_n)$ depends only on
the compatibility class of the exhaustion.
From now on, if there is no doubt as to which compatibility class
we mean, then we simply write $A_*(\Xs)$ instead of
$A_*((\Xs_n\subset\Xs)_n)$.
\elem

\blem
\label{A^* exhaustion}
Let $\Xs$ be a smooth exhaustible
algebraic stack of constant dimension $d$. 
Let $A^i(\Xs):=A_{d-i}(\Xs)$. Then $A^*:=\Oplus_iA^i(\Xs)$ has
a natural structure of a graded ring.
\elem

\bdefin
\label{f_* exhaust}
Let $\Xs$, $\Ys$ be two exhaustible algebraic stacks 
and let $(\Xs_n\subset\Xs)_n$, $(\Ys_n\subset\Ys)_n$ be 
exhaustions of $\Xs$ and $\Ys$ respectively.
Let $f:\Xs\to\Ys$ be a 
morphism which is compatible with the given
exhaustions.
Let $f_n:f^\mo(\Ys_n)\to\Ys_n$ be the restriction of $f$.
\balph
\item
Assume that $f$ is flat of relative dimension $r$.
We define the pull back morphism
$$
f^*:=\varprojlim_nf_n^*:
A_i(\Ys)
\to 
A_{i+r}(\Xs)
\quad.
$$
\item
Assume that $f$ is proper.
We define the push forward morphism
$$
f_*:=\varprojlim_n(f_n)_*:
A_i(\Xs)
\to A_i(\Ys)
\quad.
$$
\ealph
\edefin

\blem
\label{lem functorial}
Let $f:\Xs_1\to\Xs_2$ and $g:\Xs_2\to\Xs_3$ be two morphisms
between exhaustible algebraic stacks, and assume $f,g$ to be
compatible with the exhaustions. 
\balph
\item
Assume that $f,g$ are flat. Then we have $(g\comp f)^*=f^*\comp g^*$.
\item
Assume that $f,g$ are proper. Then we have $(g\comp f)_*=g_*\comp f_*$.
\item
If $f$ is an isomorphism then $f_*$ and $f^*$ are inverse to each
other.
\item
In the situation of Lemma \ref{exhaust Cart} assume that $f$ is
flat and $g$ is proper. Then we have $f^*\comp g_*=g'_*\comp(f')^*$.
\ealph
\elem

\section{The intersection ring of $\Ws$}
\label{section ring of W}

Let $\Pc$ be the category, whose objects are the finite ordered sets
\[
[1,n]:=\{1,\dots,n\}
\qquad
\text{($n$ a non-negative integer)}
\]
and whose arrows are strictly increasing mappings
between these sets.

There is a natural covariant 
functor from the category $\Pc$ to the category
of rings which to $[1,n]\in\ob(\Pc)$ associates the ring of 
Polynomials $\Q[t_1,\dots,t_n]$ and to a morphism 
$\alpha:[1,n]\to[1,m]$
in $\Pc$ associates the injective ring-homomorphism
\[
\alpha_*:
\left\{
\begin{array}{ll}
\Q[t_1,\dots,t_n] & \to  \Q[t_1,\dots,t_m] \\
t_i & \mapsto t_{\alpha(i)}
\end{array}
\right.
\]

\begin{definition}
\label{R_d}
Let $d$ be a non-negative integer.
We define $R_d$ to be the subset of the ring $\Q[t_1,\dots,t_d]$
consisting of all elements of the form
\[
\sum_{n=0}^d\left(\sum_{\alpha\in\Hom_\Pc([1,n],[1,d])}\alpha_*f_n\right)
\quad,
\]
where
\[
f_n\in\Ic_n:=\left(\prod_{i=1}^nt_i\right)\cdot\Q[t_1,\dots,t_n]
\quad\qquad
(0\leq n\leq d)
\quad.
\]
\end{definition}

\begin{lemma}
\label{lem0}
a)
The subset $R_d$ is in fact a sub-ring of $\Q[t_1,\dots,t_d]$.

b)
The kernel of the composed morphism
\[
R_d\injto\Q[t_1,\dots,t_d]\to\Q[t_1,\dots,t_d]/t_d=\Q[t_1,\dots,t_{d-1}]
\]
is $\Ic_d$ and its image is $R_{d-1}$. In particular, we have 
$R_d/\Ic_d\isomorph R_{d-1}$.
\end{lemma}

\begin{proof}
The verification is an easy exercise which we leave to the reader.
\end{proof}

We now state the main result of this section:

\begin{theorem}
\label{thm1}
The pull back homomorphism associated to the \'etale morphism
\[
[\A^1/\Gm]^d\to\Ws_d
\]
from \ref{etale} is an injective ring homomorphism
\[
A^*\Ws_d\to A^*[\A^1/\Gm]^d=\Q[t_1,\dots,t_d]
\]
whose image is the ring $R_d$ defined in \ref{R_d}.
Furthermore for every $d\geq 1$ the following diagram commutes:
$$
\xymatrix@R=2ex{
\text{$A^*\Ws_d$} \ar[d]^{\isomorph}\ar[r] 
&
\text{$A^*\Ws_{d-1}$} \ar[d]_{\isomorph} 
\\
R_d \ar[r]
&
R_{d-1}
}
$$
where the upper horizontal arrow is the pull back morphism associated to
the inclusion $\Ws_{d-1}\injto\Ws_d$ and the lower horizontal arrow
is the one defined in Lemma \ref{lem0}.
\end{theorem}

The Theorem will follow from Lemma \ref{lem1} through \ref{lem5}
below.

\begin{lemma}
\label{lem1}
There is a natural commutative diagram of stacks,
\[
\xymatrix{
\text{$\BGm^d$} \ar@{^(->}[r]^(.4){i_1} \ar@{=}[d] & 
\text{$[\A^1/\Gm]^d$} \ar[d]^{f_d}  &
\text{$[\A^d\setminus\{0\}/\Gm^d]$} 
\ar@{_(->}[l]_(.55){j_1} \ar[d]^{g_{d-1}}\\
\text{$\BGm^d$} \ar@{^(->}[r]^{i_2} & 
\text{$\Ws_d$} &
\text{$\Ws_{d-1}$} \ar@{_(->}[l]_{j_2}
}
\]
with the following properties
\balph
\item
For $\nu=1,2$ the morphism $i_\nu$ is a closed immersion,
and $j_\nu$ is an open immersion whose image is complementary
to the image of $i_\nu$.
\item
The morphism $f_d$ is the \'etale morphism from 
Proposition \ref{etale}.
\item
The left square is Cartesian.
\ealph
\end{lemma}

\bpf
The morphisms $i_1$, $j_1$ are induced from the
following morphisms of schemes
$$
\xymatrix{
\Spec(k)=\{0\} \ar@{^(->}[r] & \text{$\A^d$} & 
 \text{$\A^d\setminus\{0\}$} \ar@{_(->}[l]
}
$$
by passing to quotients.
Let $T$ be a $k$-scheme. On the level of $T$ valued points
the morphisms $i_2$, $j_2$ are defined
as follows:
\beq
i_2 &:&(L_1,\dots,L_d) \mapsto (L_1,0)\dots(L_d,0)\quad,\\
j_2 &:&(L_1,\lambda_1)\dots(L_{d-1},\lambda_{d-1})\mapsto
(L_1,\lambda_1)\dots(L_{d-1},\lambda_{d-1})\quad.
\eeq
Finally, we define the morphism $g_{d-1}$.
A $T$ valued point of 
$[\A^d\setminus\{0\}/\Gm^d]$
is a tuple $\xi=((L_1,\lambda_1),\dots,(L_d,\lambda_d))$ of $d$
line bundles with sections, such that the open sets
$T_i:=\{\lambda_i\neq 0\}$ form a covering of $T$.
For each $i$ let $\eta_i\in\Ws_{d-1}(T_i)$ be the object which
arises from $\xi|_{T_i}$ by cancelling $(L_i,\lambda_i)|_{T_i}$
and all the other $(L_j,\lambda_j)|_{T_i}$ for 
which $\lambda_j|_{T_i}$ is invertible.
Obviously the $\eta_i$ glue together to form an object 
$\eta\in\Ws_{d-1}(T)$ and we let $g_{d-1}(\xi):=\eta$.

With these definitions it is immediate that the diagram
has the stated properties.
\epf

\blem
\label{lem2}
The diagram in Lemma \ref{lem1} induces the following 
commutative diagram
with exact rows:
$$
\xymatrix@C=7ex{
0 \ar[r] 
&
\Q[t_1,\dots,t_d] \ar@{=}[d] \ar[r]^{\prod_it_i} 
& 
\Q[t_1,\dots,t_d] \ar@{=}[d] \ar[r]
& 
\Q[t_1,\dots,t_d]/(\prod_i t_i) \ar@{=}[d] \ar[r]
&
0
\\
0 \ar[r] 
&
\text{$A^*(\BGm^d)$} \ar@{=}[d] \ar[r]^{(i_1)_*}  
&
\text{$A^*([\A^1/\Gm]^d)$} \ar[r]^{j_1^*}
& 
\text{$A^*([\A^d\setminus\{0\}/\Gm^d])$} \ar[r]
&
0
\\
0 \ar[r] 
&
\text{$A^*(\BGm^d)$}  \ar[r]^{(i_2)_*}  
&
\text{$A^*(\Ws_d)$} \ar[r]^{j_2^*} \ar[u]^{f_d^*} 
&
\text{$A^*(\Ws_{d-1})$} \ar[r] \ar[u]^{g_{d-1}^*}
&
0
}
$$
\elem

\bpf
For $i=1,\dots,d$ let $V_i$ be a $k$-vector space of dimension $p$
on which $\Gm$ operates by scalar multiplication, and let
$U_i:= V_i\setminus\{0\}$. Then $\Gm^d$ acts freely
on $U:=U_1\times\dots\times U_d$ and the complement of $U$ in
$V:=V_1\oplus\dots\oplus V_d$ is of codimension $p$.
Thus we have (cf. \cite{Totaro}, \S 1)
$$
A^i(\BGm^d)=A^{i}(U/\Gm^d)=A^{i}(\Pp)
=(\Q[t_1,\dots,t_d]/(t_1^p,\dots,t_d^p))_{i}
=(\Q[t_1,\dots,t_d])_{i}
$$
for $i<p$, where 
$\Pp:=\Pp(V_1)\times\dots\times\Pp(V_d)$
and $t_i\in A^{1}(\Pp)$
denotes the pull-back of the class of a hyperplane in $\Pp(V_i)$.

The stack $\BGm^d$ classifies split vector bundles of rank $d$,
$q:[\A^1/\Gm]^d\to \BGm^d$ being the universal bundle and
$i_1:\BGm^d\to[\A^1/\Gm]^d$ being the zero section.
A finite dimensional approximation of $q$ and $i_1$ is given
by 
$$
\pi:E\to \Pp
\quad\text{and}\quad
s:\Pp\to E
$$
where $E=E_1\boxtimes\dots\boxtimes E_d$ and $E_i=\Oo_{\Pp(V_i)}(1)$
is the dual of the tautological line bundle on $\Pp(V_i)$.
By \cite{Fulton} Theorem 3.3 and Example 3.3.2,
$\pi^*$ identifies $A^*(E)$ with
$A^*(\Pp)=\Q[t_1,\dots,t_d]/(t_1^p,\dots,t_d^p)$
and $s_*:A^*(\Pp)\to A^{*+d}(E)$ translates into multiplication
with the element $c_d(E)=\prod_it_i$.
Letting $p$ grow to infinity, we obtain the commutative square
in the left upper corner of the diagram.

The two lower rows are excision sequences and thus exact to the
right (\cite{Kresch}, Theorem 2.1.12 (iv)).
Since the left square in Lemma \ref{lem1} is Cartesian, 
we have $(i_1)_*=f_d^*\comp(i_2)_*$. By our explicit description
above the morphism $(i_1)_*$ is injective. Therefore also
$(i_2)_*$ is injective.
\epf

\bdefin
\label{smaller}
\balph
\item
For a strictly increasing map $\alpha:[1,n]\to[1,m]$
we let
$$
\alpha^*:\Q[t_1,\dots,t_m]\to\Q[t_1,\dots,t_n]
$$
be the ring homomorphism defined by
$$
\alpha^*(t_i)=
\bca
t_{j} &\text{if there is a $j\in[1,n]$ with $\alpha(j)=i$,}\\
0 & \text{else.}
\eca
$$
By abuse of notation for $I\subset\ker\alpha^*$
we will denote by the same symbol $\alpha^*$
the induced morphism $\Q[t_1,\dots,t_m]/I\to\Q[t_1,\dots,t_n]$.
\item
Let $\alpha_1,\alpha_2:[1,n]\injto[1,m]$ be two strictly increasing
maps. We will say that $\alpha_1$ is {\em smaller} (resp. {\em larger} 
than 
$\alpha_2$, if the sum over all values of $\alpha_1$ is
smaller (resp. larger) than the sum over all values of $\alpha_2$.
\ealph
\edefin

\blem
\label{lem3}
\broman
\item
For $i\in[1,d+1]$ let 
$h_{d,i}$
be the composite morphism
$$
[\A^1/\Gm]^{d}=[(\A^{i-1}\times\Gm\times\A^{d-i+1})/\Gm^{d+1}]
\injto
[(\A^{d+1}\setminus\{0\})/\Gm^{d+1}]
$$
Then the following diagram of stack morphisms commutes:
$$
\xymatrix{
&
\text{$[\A^1/\Gm]^{d}$} \ar[dr]^{f_{d}} \ar[dl]_{h_{d,i}}
&
\\
\text{$[(\A^{d+1}\setminus\{0\})/\Gm^{d+1}]$} \ar[rr]^(0.6){g_{d}}  
&
&
\text{$\Ws_{d}$} 
}
$$

\item
Let $\alpha:[1,d]\injto[1,d+1]$ be the strictly
increasing map, which omits $i\in[1,d+1]$.
Then the following diagram of ring homomorphisms commutes:
$$
\xymatrix{
\Q[t_1,\dots,t_{d+1}]/\prod_{i=1}^{d+1}t_i 
\ar@{=}[d] \ar[r]^(.6){\alpha^*}
&
\Q[t_1,\dots,t_d] \ar@{=}[d]
\\
\text{$A^*([(\A^{d+1}\setminus\{0\})/\Gm^{d+1}])$} 
\ar[r]^(.65){h_{d,i}^*}
&
\text{$[\A^1/\Gm]^{d}$}
}
$$
\item
For $d\geq 1$ the morphisms $f_d^*$ and $g_{d-1}^*$ in Lemma \ref{lem2}
are injective.
\eroman
\elem

\bpf
(i):
Let $T$ be a $k$-scheme.
On $T$-valued points the morphism $h_{d,i}$ is
defined by 
$$
((L_1,\lambda_1),\dots,(L_d,\lambda_d))\mapsto
((L_1,\lambda_1),\dots,(L_{i-1},\lambda_{i-1}),(\Oo_T,1),
(L_{i+1},\lambda_{i+1}),\dots,(L_d,\lambda_d))
\quad.
$$
From this the first part of the lemma follows immediately.

(ii):
This is clear since the morphism 
$$
h_{d,i}':
[\A^1/\Gm]^{d}=[(\A^{i-1}\times\Gm\times\A^{d-i+1})/\Gm^{d+1}]
\injto
[\A^{d+1}/\Gm^{d+1}]
$$
factorises through $h_{d,i}$,
and we obviously have $(h_{d,i}')^*=\alpha^*$.

(iii):
This follows by induction on $d$: The morphisms
$f_1$ and $g_0$ and therefore also $f^*_1$ and $g^*_0$ 
are isomorphisms.
Assume that injectivity of $f^*_d$ and $g^*_{d-1}$
holds for some $d\geq 1$.
By part (i) we have $f^*_d=(h_{d,i})^*\comp g^*_d$,
thus injectivity of $g^*_d$ follows.
But from the diagram in Lemma \ref{lem1} it is clear
that injectivity of $g^*_d$ implies injectivity of
$f^*_{d+1}$.
\epf

\blem
\label{lem4}
For any ring $A$ the following sequence (of sets)
is exact:
$$
\def\objectstyle{\displaystyle}
\xymatrix{
*{\Hom\left(A,\Q[t_1,\dots,t_d]/t_1\dots t_d\right)} 
\ar[r]^(.45)u
&
*{\prod_{\alpha:[1,d-1]\injto[1,d]}\Hom(A,\Q[t_1,\dots,t_{d-1}])}
\\
\ar@<1.5ex>[r]^(.3){v_1} \ar@<-.5ex>[r]^(.3){v_2}
&
*{\prod_{\beta:[1,d-2]\injto[1,d]}\Hom(A,\Q[t_1,\dots,t_{d-2}])}
}
$$ 
where the mappings $u,v_1,v_2$ are defined as follows:
\begin{itemize}
\item
The $\alpha$-th component of
$u(\varphi)$ is $\alpha^*\comp\varphi$.
\item
The $\beta$-th component of
$v_1((\varphi_\alpha)_\alpha)$
(resp. $v_2((\varphi_\alpha)_\alpha)$) is
$\gamma^*\comp\varphi_\alpha$, where
$\alpha$ is the smaller (resp. larger, 
cf. Definition \ref{smaller})
of the two strictly increasing maps $[1,d-1]\injto[1,d]$ such that there
exists a $\gamma:[1,d-2]\injto[1,d-1]$ with
$\beta=\alpha\comp\gamma$.
\end{itemize}
\elem

\bpf
We leave the easy verification to the reader.
Note that the lemma means geometrically that a morphism to some scheme 
from a union of 
hyperplanes in an affine space is the same thing as
a tuple consisting of morphisms to the scheme
from each hyperplane which coincide on the intersections. 
\epf

\blem
\label{lem5}
The \'etale morphism $[\A^1/\Gm]^d\to\Ws_d$
induces an injective ring homomorphism 
\[
A^*\Ws_d\to A^*[\A^1/\Gm]^d=\Q[t_1,\dots,t_d]
\]
whose image is the ring $R_d$ defined in \ref{R_d}.
\elem

\bpf
From Lemma \ref{lem3} we already know that 
$f_d^*:A^*(\Ws_d)\to\Q[t_1,\dots,t_d]$ is injective
for all $d\geq 1$.
Therefore we can identify $A^*(\Ws_d)$ with
a sub-ring of $\Q[t_1,\dots,t_d]$.
Since $f_0$, $f_1$ are isomorphisms, and 
$R_0=\Q$, $R_1=\Q[t_1]$ it follows that
$A^*(\Ws_d)=R_d$ for $d\leq 1$.

We proceed by induction over $d$.
Assume that $A^*(\Ws_d)=R_d$ for some $d\geq 1$.
By Lemma \ref{lem3} the homomorphism
$$
g^*_d:R_d=A^*(\Ws_d)\to\Q[t_1,\dots,t_{d+1}]/t_1\dots t_{d+1}
$$ 
has the property that
for all strictly increasing $\alpha:[1,d]\injto[1,d+1]$
we have
$\alpha^*\comp g^*_d=f_d^*$
and Lemma \ref{lem4} implies that $g^*_d$ is uniquely 
determined by this property.
We claim that for any $n\in[0,d]$ and 
$f\in\left(\prod_{i=1}^nt_i\right)\Q[t_1,\dots,t_n]$
we have
$$
g^*_d\left(\sum_{\gamma:[1,n]\injto[1,d]}\gamma_*(f)\right)=
\sum_{\gamma:[1,n]\injto[1,d+1]}\gamma_*(f)
\quad\mod \prod_{i=1}^{d+1}t_i
\quad.
$$
Indeed, this follows from the above uniqueness statement for $g^*_d$
and the fact that for all
strictly increasing $\alpha:[1,d]\injto[1,d+1]$
we have
$$
\alpha^*\left(\sum_{\gamma:[1,n]\injto[1,d+1]}\gamma_*(f)\right)
=
\sum_{\gamma:[1,n]\injto[1,d]}\gamma_*(f)
\quad.
$$
This shows in particular that 
the image of $g^*_d$ in $\Q[t_1,\dots,t_{d+1}]/(t_1\dots t_{d+1})$
is $R_{d+1}/\I_{d+1}$.
By the diagram in Lemma \ref{lem2} it follows that
the image of $f^*_{d+1}$ in $\Q[t_1,\dots,t_{d+1}]$ 
is precisely the ring $R_{d+1}$.
This proves the Lemma.
\epf

This finishes the proof of Theorem \ref{thm1} and thus of 
the isomorphism $A^*(\Ws_d)\isomorph R_d$ of graded rings.
Next we compute $A^*(\Ws)$. Note that $\Ws$ is not of finite
type, but it is exhaustible in the sense of \S\ref{section exhaustible}
since $(\Ws_n\subset\Ws)_n$ is an exhaustion for $\Ws$.
Therefore $A^*(\Ws)$ in the sense of Definition \ref{A_* exhaust}
is defined: $A^*(\Ws)=\varprojlim'_nA^*(\Ws_n)=\Oplus_i\varprojlim_n(A^i(\Ws_n))$.

We will prove below that $A^*(\Ws)$ is isomorphic to a ring
$R$ which is defined as follows.
\label{R}
Recall from Definition \ref{R_d} the definition of the ideal
$\Ic_n$ and as a group let $R:=\bigoplus_{n\geq 0}\Ic_n$.
We define a multiplication map $R\tensor R\to R$ as follows:
The product of two elements
$f=(f_n)_{n\geq 0}$
and
$g=(g_n)_{n\geq 0}$
in $R$ is the element $h=(h_n)_{n\geq 0}$
where
\[
h_n:=\sum_{\alpha,\beta}
(\alpha_*f_{p})(\beta_*g_{q})
\quad.
\]
In this sum $(\alpha,\beta)$ runs through all elements
$\Hom_\Pc([1,p],[1,n])\times\Hom_\Pc([1,q],[1,n])$
such that the (not necessarily disjoint)
union of the images of $\alpha$ and $\beta$
is the whole of $[1,n]$.
Let the grading of $R$ be induced by the grading
of the $\Ic_n\subset\Q[t_1,\dots,t_n]$.
It is easy to check that with 
these structures $R$ is a commutative associative
graded ring with unity. 

\bthm
\label{thm2}
There is a natural isomorphism $A^*(\Ws)\isomto R$ of graded rings.
\ethm

\bpf
For $d\geq 1$ there is a natural morphism of graded rings
$R\to R_d$ which maps
an element $(f_n)_{n\geq 0}$ to 
$$
F:=
\sum_{n=0}^d\left(
\sum_{\alpha\in\Hom_\Pc([1,n],[1,d])}\alpha_*f_n\right)
\quad.
$$
These morphisms induce a morphism $R\to\varprojlim'_nR_n$
where the inverse limit is taken in the category of graded rings 
with respect to the morphisms $R_n\to R_{n-1}$ from \ref{lem0}.
By Theorem \ref{thm1} it suffices to show that this is an isomorphism.
The main point in the verification of this is that
a polynomial $F\in R_d$ as above uniquely determines the polynomials
$f_n$ for $0\leq n\leq d$.
Indeed, if we write $F=\sum_{I\in\N_0^d}a_It^I$ as sum of monimials
(in multi-index notation), then
$$
f_n=\sum_Ia_It^I
\quad,
$$
where $I$ runs through all tuples $(i_1,\dots,i_d)$
with $i_1,\dots,i_n>0$ and $i_{n+1}=\dots=i_d=0$.
We leave the details to the reader.
\epf

\bco
\label{A^*(XxW)}
Let $\Xs$ be a smooth exhaustible stack (cf. \S\ref{section exhaustible})
with exhaustion 
$(\Xs_n\subset\Xs)_n$.
Let the stack $\Xs\times \Ws$ be endowed with the exhaustion 
$(\Xs_n\times\Ws_n\subset\Xs\times\Ws)_n$.
Then there is a natural
isomorphism $A^*(\Xs\times\Ws)\isomto A^*(\Xs)\tensor R$ of graded rings.
\eco

\bpf
Note first of all that for any $n$ there is a natural isomorphism 
$A^*(\Xs_n\times[\A^1/\Gm])\isomto A^*(\Xs_n)\tensor\Q[t]$.
This follows essentially from the following well known fact:
Let $N>0$ and let $E_N\to\Pp^N$ be the tautologycal line bundle 
on the $N$-dimensional projective space. Then there is a natural
isomorphism $A^*(\Xs_n\times E_N)\isomto A^*(\Xs_n)\tensor\Q[t]/(t^N)$.

It follows from this observation that we have 
$A^*(\Xs_n\times\Ws_n)=A^*(\Xs_n)\tensor R_n$
(simply repeat the proof of Lemmas \ref{lem1}-\ref{lem5} with
$\Xs_n\times \Ws_n$ in place of $\Ws_n$, with $\Xs_n\times \BGm$
in place of $\BGm$, with $\Xs_n\times[A^1/\Gm]^n$ in place of $[A^1/\Gm]^n$,
etc.

Let $R^{(q)}_n\subset R_n$ and $R^{(q)}\subset R$
be the parts of degree $q$. Then we have $R^{(q)}=R^{(q)}_n$
for $n>q$ and it follows
\begin{multline*}
A^i(\Xs\times \Ws)=
\varprojlim_nA^i(\Xs_n\times\Ws_n)=
\varprojlim_n\bigoplus_{p+q=i}A^p(\Xs_n)\tensor R^{(q)}_n
=\\=
\bigoplus_{p+q=i}\varprojlim_n(A^p(\Xs_n))\tensor R^{(q)}=
\bigoplus_{p+q=i}A^p(\Xs)\tensor R^{(q)}
\end{multline*}

\epf

\section{Zollstocks and the Artin stack $\Moz$}
\label{section zollstock}

\bdefin
\balph
\item
A {\em zollstock} is a two-pointed nodal curve $(C,x_1,x_2)$
of genus zero with irreducible components $R_1,\dots,R_l$
such that
\begin{itemize}
\item
$x_1\in R_1$, $x_2\in R_l$,
\item
for $i=1,\dots,l-1$
the intersection
$R_i\cap R_{i+1}$ consists of one point $y_i$,
\item 
all the other intersections between the $R_1,\dots,R_l$
are empty,
\item
the points $x_1,y_1,\dots,y_{l-1},x_2$ are distinct.
\end{itemize}
The number $l$ of its components is called the {\em length} of the
zollstock.
\item
Let $S$ be scheme. A {\em zollstock  over $S$} is a 
diagram
$
\xymatrix{
X \ar[r]_\pi
&
S \ar@<-.3ex>@/_/[l]_{s_1,s_2}
}
$
where $\pi$ is a proper
flat morphism
and $s_1,s_2$ are two sections of $\pi$, such
such that all fibres of $\pi$ are zollstocks
(of possibly variable length).
\ealph
\edefin

``Zollstock'' is the German word for ``folding rule''.
We have chosen this name, because the usual schematic
picture for curves matching the above definition 
reminds one of this device:
$$
\vcenter{
\xy <1cm,0cm>:
(0.8,0.5);
(2.2,1) **@{-};
?<<< *{\bullet}
;
(1.8,1);
(3.2,0.5) **@{-};
;
(2.8,0.5);
(4.2,1) **@{-};
;
(4.5,0.5)*{\cdot};
(4.8,0.5)*{\cdot};
(5.1,0.5)*{\cdot};
(5.4,0.5)*{\cdot};
;
(5.8,1);
(7.2,0.5) **@{-};
;
(6.8,0.5);
(8.2,1) **@{-}
;
(7.8,1);
(9.2,0.5) **@{-}
?>>> *{\bullet}
\endxy
}
$$

\bex
\label{zollstock ex}
\barabic
\item
Let $L_1$ and $L_2$ be two line bundles over $S$.
Let $\Cc:=\Pp(L_1\oplus L_2)$ and let
$s_i:S\to \Cc$ be the morphisms defined by the invertible quotient
$L_1\oplus L_2\to L_i$ for $i=1,2$.
Then $(\Cc,s_1,s_2)$ is a zollstock of constant length $1$
over $S$, which we denote by the symbol $|L_1,L_2|$.
For later reference we remark that we have
$s_1^*\Oo_\Cc(s_1)=L_1\tensor L_2^\mo$
and
$s_2^*\Oo_\Cc(s_2)=L_1^\mo\tensor L_2$.

\item
Given a zollstock 
$(X\to S, s_1,s_2)$ over $S$ 
and a Divisor $D\subset S$ we obtain a new zollstock
$(X'\to S, s'_1,s'_2)$
by blowing up $X$ in $s_1(D)$ and taking $s'_i$ to be
the proper transform of $s_i$.

\item
More generally, if $(X\to S, s_1,s_2)$ is a zollstock over $S$ 
and $(X', x')\to(X,s_2)$ is a modification of $X$ at $s_2$
then $(X'\to S, s_1, x')$ is also a zollstock over $S$.
If $(X', x')\to(X,s_2)$ is a simple modification corresponding
via Proposition \ref{T1 isom A1/Gm} to a line bundle with
section $(L,\lambda)$, then we write
$
(s_1,X,s_2)\vdash(L,\lambda)
$
for this zollstock. 
The zollstock 
$
(L,\lambda)\dashv(s_1,X,s_2)
$
is defined analogously, interchanging the role of $s_1$ and $s_2$.

\item
If $(X\to S, s_1,s_2)$
and $(Y\to S, t_1,t_2)$
are two zollstocks then $X\bot Y:=X\sqcup Y/(s_2=t_1)$ is
a zollstock (over $S$).

\item
Escher's staircase $X\to S$ on page \pageref{escher} is a zollstock over $S$.
\earabic
\eex

For the convenience of the reader we recall the following 
useful result
from \cite{factorization}, \S2.

\blem
\label{usefull}
Let $S$ be a scheme, let $L_0,\dots,L_k$, $M_0,\dots,M_l$ be
line bundles over $S$ and let $\lambda_i\in H^0(S,L_i)$,
$\mu_j\in H^0(S,M_j)$ be sections for $1\leq i\leq k$, $1\leq j\leq l$.
Then there are canonical isomorphisms of zollstocks over $S$ as follows:
\balph
\item
\begin{multline*}
\Bigl(
|L_0,\Oo_X|\vdash(L_1,\lambda_1)\vdash\dots\vdash(L_k,\lambda_k)
\Bigr)
\bot
\Bigl(
|M_0,\Oo_X|\vdash(M_1,\mu_1)\vdash\dots\vdash(M_l,\mu_l)
\Bigr)
\isomorph\\\isomorph
|L_0,\Oo_X|\vdash(L_1,\lambda_1)\vdash\dots\vdash(L_k,\lambda_k)
\vdash
\left(M_0\tensor\Tensor_{i=0}^kL_i^\mo,0\right)
\vdash(M_1,\mu_1)\vdash\dots\vdash(M_l,\mu_l)
\end{multline*}
\item
$$
|L_0,\Oo_X|\vdash(L_1,\lambda_1)\vdash\dots\vdash(L_k,\lambda_k)
\isomorph
(L_1,\lambda_1)\dashv\dots\dashv(L_k,\lambda_k)
\dashv
|\Oo_S,\Tensor_{i=0}^kL_i^\mo|
$$
\ealph
\elem

\bpf
The first isomorphism follows easily from Lemma 2.9 in \cite{factorization}.
The second is shown in Lemma 2.11 in loc. cit.
\epf

Let $\Zs$ be the $\C$-groupoid parametrising zollstocks.
For $l\geq 1$ let $\Zs_l$ be the subgroupoid of $\Zs$
parametrising zollstocks of length at most $l$.

\bpr
\label{M02}
There is a natural isomorphism
$$
\Moz:=\BGm\times\Ws\to\Zs
$$
of $\C$-groupoids.
It maps the subgroupoid $\BGm\times\Ws_l$ isomorphically to the
subgroupoid $\Zs_{l+1}$.
\epr

\bpf
For a scheme $S$ let
$(\BGm\times\Ws)(S)\to\Zs(S)$
be defined as follows:
An object in $(\BGm\times\Ws)(S)$ is a pair $(L,\xi)$ consisting
of a line bundle $L$ over $S$ and and object $\xi\in\Ws(S)$.
Let $(Y:=\Pp(L\oplus\Oo_S)\to S   ,s_0,s_\infty)$ be the zollstock over $S$
as in example \ref{zollstock ex}(1) and let $(X,x_2)\to(Y,\infty)$
be the modification of $Y$ at $\infty$ defined by $\xi$ via
the isomorphism $\Ws(S)\to\Tails(Y,s_\infty)(S)$
in Proposition \ref{BxW isom T} . 
Then $(X\to S,x_1=s_0,x_2)$ is a zollstock over $S$, i.e. an object
in $\Zs(S)$, which we take to be the image of $(L,\xi)$.

In the reverse direction let $(\pi:X\to S, x_1, x_2)$ be a zollstock over
$S$. Then $E:=\pi_*\Oo_X(x_1)$ is a locally free sheaf of rank two over
$S$ and the natural morphism $\pi^*E\to\Oo_X(x_1)$ is surjective - thus
defining an $S$-morphism $f:X\to\Pp(E)$. Let $s_i:=f\comp x_i$. Then
it is easy to see that $(X,x_2)\to(\Pp(E),s_2)$ is a modification of $\Pp(E)$
at $s_2$. Let $\xi\in\Ws(S)$ be the object which corresponds to this
modification by Proposition \ref{BxW isom T} and let $L:=x_1^*\Oo_X(x_1)$.
We define $\Zs(S)\to(\BGm\times\Ws)(S)$ to be the morphism which 
to $(X\to S, x_1, x_2)$ associates the object $(L,\xi)$.

It is not difficult to see that both constructions commute with base
change and define morphisms of $\C$-groupoids which are quasi-inverse
to each other. The last part of the proposition is immediate from
the constructions.
\epf

By Proposition \ref{M02} and Corollary \ref{etale cor} we have:
\bco
\label{M02 cor}
The groupoids $\Zs$ and $\Zs_l$ are smooth algebraic stacks
of dimension $-1$. Moreover the $\Zs_l$ are of finite type
and we have the following inclusion of open substacks
$$
\BGm=\Zs_1\subset\Zs_2\subset\dots\subset\Zs_l\subset\dots\subset\Zs
\quad.
$$
\eco

We have natural morphisms
\begin{align*}
\clutch:&\ \Zs\times\Zs  \to\Zs\quad,\\
\sigma:&\ \Zs  \isomto\Zs\quad.
\end{align*}
The morphism $\clutch$ is the clutching morphism, which to 
two zollstocks $(X,x_1,x_2)$ and $(Y,y_1,y_2)$ associates the
zollstock $(X\sqcup Y/(x_2=y_1),x_1,y_2)$ which arises by identifying
the second section of the first zollstock with the first section of
the second zollstock.
The other morphism $\sigma$ is the obvious involution 
which comes from interchanging the two sections $x_1,x_2$ in a zollstock
$(X,x_1,x_2)$.
The following Lemma is immediate from the definitions:

\blem
\label{Z monoid}
The clutching morphism  provides $\Zs$ with the structure of a
monoid in the category of stacks. 
Furthermore, the isomorphism $\sigma:\Zs\to\Zs$
is an anti-isomorphism of that monoid.
More precisely the 
following diagrams are canonically $2$-commutative:
$$
\xymatrix{
\text{$\Zs\times\Zs\times\Zs$} 
\ar[r]^(0.55){\id\times\clutch}
\ar[d]_{\clutch\times\id}
&
\text{$\Zs\times\Zs$}
\ar[d]^{\clutch}
\\
\text{$\Zs\times\Zs$}
\ar[r]^{\clutch}
&
\text{$\Zs$}
}
\qquad
\xymatrix{
\text{$\Zs\times\Zs$}
\ar[rr]^{\clutch}
\ar[d]^{\exchange}
& &
\text{$\Zs$}
\ar[d]^{\sigma}
\\
\text{$\Zs\times\Zs$}
\ar[r]^{\sigma\times\sigma}
&
\text{$\Zs\times\Zs$}
\ar[r]^\clutch
&
\text{$\Zs$}
}
$$
\elem

Via the isomorphism in Proposition \ref{M02} the clutching morphism
$\clutch$ and the involution $\sigma$ induce
morphisms
\label{clutch sigma}
\begin{align*}
\clutch:&\ \Moz\times\Moz  \to\Moz \quad,\\
\sigma:&\ \Moz  \isomto\Moz \quad,
\end{align*}
which we denote by the same symbol.
Explicitly these arrows are given by
\begin{multline*}
\clutch:
\Bigl(\bigl(L_0,(L_1,\lambda_1)\dots(L_k,\lambda_k)\bigr),
\bigl(M_0,(M_1,\mu_1)\dots(M_l,\mu_l)\bigr)\Bigr)
\\
\mapsto
\left(L_0,(L_1,\lambda_1)\dots(L_k,\lambda_k)
\left(M_0\tensor\Tensor_{i=0}^kL_i^{-1},0\right)
(M_1,\mu_1)\dots(M_l,\mu_l)\right)
\end{multline*}
and
$$
\label{sigma}
\sigma:\Bigl(L_0,(L_1,\lambda_1)\dots(L_k,\lambda_k)\Bigr)
\mapsto
\left(\Tensor_{i=0}^kL_i^{-1},(L_k,\lambda_k)\dots(L_1,\lambda_1)\right)
\quad.
$$
This follows from Lemma \ref{usefull}.
Of course the morphisms $\clutch$ and $\sigma$ restrict to the open
substacks $\Zs_l$  parameterising zollstocks of
bounded length:
\begin{align*}
\clutch_{k,l}&:\Zs_k\times\Zs_l\to\Zs_{k+l}
\quad,
\\
\sigma_l&:\Zs_l\isomto\Zs_l
\quad.
\end{align*}

The clutching morphism can be looked upon from yet another
point of view. To explain this we need to introduce a notation.
For a flat morphism $\pi:X\to Y$ of schemes (or algebraic stacks) we let
$\Sing(X/Y)\subset X$ be the closed subscheme (or substack) where $\pi$ fails to
be smooth. It is   
defined by the first Fitting ideal of the sheaf $\Omega^1_{X/Y}$.

\bpr
\label{sing}
Let $(\Xs\to\Zs,\x_1,\x_2)$ be the universal zollstock over $\Zs$. Then there 
is a canonical isomorphism $\Sing(\Xs/\Zs)\isomto\Zs\times\Zs$
such that the following diagram commutes:
$$
\xymatrix{
\text{$\Sing(\Xs/\Zs)$}\ar[d]_\isomorph\ar@{^(->}[r]
&
\text{$\Xs$} \ar[d]
\\
\text{$\Zs\times\Zs$} \ar[r]^\clutch
&
\text{$\Zs$}
}
$$
\epr

\bpf
Let $S$ be a scheme and let 
$\xi:S\to\Zs$ be given by a zollstock $(\pi:X\to S,x_1,x_2)$ over $S$.
Let $S':=\Sing(X/S)$ and let $X':=X\times_SS'$.
Let $\nu:X'_1\sqcup X'_2\to X'$ be the partial normalisation of $X'$
along the canonical section $y':S'\to X'$ of $X'\to S'$.
There are natural
sections $x'_1,y'_1:S'\to X'_1$ and $x'_2,y'_2:S'\to X'_2$
such that $x_1|_{S'}=\nu\comp x'_1$, $x_2|_{S'}=\nu\comp x'_2$, 
$y'=\nu\comp y'_1=\nu\comp y'_2$. Thus we have zollstocks
$\xi'_1=(X'_1\to S',x'_1,y'_1)$ and
$\xi'_2=(X'_2\to S',y'_2,x'_2)$ over $S'$ such that 
$\clutch(\xi'_1,\xi'_2)=\xi|_{S'}$.
The proposition follows from the easily verifiable fact that the
diagram
$$
\xymatrix{
\text{$S'$}\ar[d]_{\pi|_{S'}}\ar[r]^(.4){(\xi'_1,\xi'_2)}
&
\text{$\Zs\times\Zs$} \ar[d]^\clutch
\\
S \ar[r]^\xi
&
\text{$\Zs$}
}
$$
is Cartesian.
\epf

For the construction of the extended operad (cf. \S\ref{section operad})
it will be important that the monoid stack $\Zs$ acts on the stack
of modifications of a pointed curve $(C/B,x)$ over a base $B$. 
Namely there is a natural morphism
$$
\label{Z acts}
\clutch_{(C/B,x)}:\Tails(C/B,x)\times\Zs\to\Tails(C/B,x)
$$
which is defined as follows.
Let $S$ be a $B$-scheme, let $(C'/S,x')\to(C_S/S,x_S)$ be an
object in $\Tails(C/B,x)(S)$ and let $(X/S,x_1,x_2)$ be a zollstock over $S$. 
The image of this pair by $\clutch_{(C/B,x)}$
is the modification $(C''/S,x'')\to(C_S/S,x_S)$ where 
$$
(C'',x''):=(X\sqcup C'/(x_2=x'),x_1)\quad,
$$
and the morphism $C''\to C_S$ is the morphism $C''\to C'$ which
contracts $X$ to $x'$, followed by the morphism $C'\to C_S$.
With this definition it is clear that the following diagram commutes:
$$
\xymatrix{
\text{$\Tails(C/B,x)\times\Zs\times\Zs$}
\ar[d]_{\id\times\clutch}
\ar[rr]^(.55){\clutch_{(C/B,x)}\times\id}
&&
\text{$\Tails(C/B,x)\times\Zs$}
\ar[d]^{\clutch_{(C/B,x)}}
\\
\text{$\Tails(C/B,x)\times\Zs$} 
\ar[rr]^{\clutch_{(C/B,x)}}
&&
\text{$\Tails(C/B,x)$}
}
$$

Via the isomorphisms of Prop. \ref{M02} and \ref{BxW isom T}
$\clutch_{(C/B,x)}$ induces a right action of the stack $\Moz$ 
on the stack $B\times\Ws$. We are interested in a concrete description
of this action.

\bdefin
\label{cK def}
For a line bundle $K$ on $B$
let $\clutch_K:(B\times\Ws)\times\Moz\to(B\times\Ws)$ be the morphism defined
by
\begin{multline*}
\Bigl((L_1,\lambda_1)\dots(L_r,\lambda_r),
\bigl(M_0,(M_1,\mu_1)\dots(M_s,\mu_s)\bigr)\Bigr)
\mapsto \\
(L_1,\lambda_1)\dots(L_r,\lambda_r)
\left(M_0\tensor K\tensor\Tensor_{i=1}^rL_i^\mo,0\right)
(M_1,\mu_1)\dots(M_s,\mu_s)
\end{multline*}
\edefin

By Lemma \ref{usefull}, the
following diagram commutes
$$
\label{cK diag}
\xymatrix{
\text{$(B\times\Ws)\times\Moz$}
\ar[d]_{\text{Prop.\ref{BxW isom T}}}^{\text{Prop.\ref{M02}}}
\ar[rr]^{\clutch_K}
&&
\text{$(B\times\Ws)$}
\ar[d]^{\text{Prop.\ref{BxW isom T}}}
\\
\text{$\Tails(C/B,x)\times\Zs$}
\ar[rr]^{\clutch_{(C/B,x)}}
&&
\text{$\Tails(C/B,x)$}
}
$$
where $K$ is the line bundle $x^*\Oo_C(x)$ on $B$.

Here is the analogue of Proposition \ref{sing}:
\bpr
\label{singt}
Let $(C/B,x)$ be a pointed curve over a scheme $B$.
For abbreviation let us write $\Tails$ for the $B$-stack $\Tails(C/B,x)$.
Let $(\Xs/\Tails,y)\to(C_\Tails/\Tails,x_\Tails)$ be the universal modification
over $\Tails$. Let $Sing'(\Xs/\Tails)$ be defined by the Cartesian
diagram
$$
\xymatrix{
\text{$\Sing'(\Xs/\Tails)$}\ar[d]\ar[r]
&
\text{$\Sing(\Xs/\Tails)$} \ar@{_(->}[d]
\\
\text{$\Tails$} \ar@{^(->}[r]^{x_\Tails}
&
\text{$\Xs$}
}
$$
Then there is
is a canonical isomorphism 
$\Sing'(\Xs/\Tails)\isomto\Tails\times\Zs$
such that the following diagram commutes:
$$
\xymatrix{
\text{$\Sing'(\Xs/\Tails)$}\ar[d]_\isomorph\ar@{^(->}[r]
&
\text{$\Xs$} \ar[d]
\\
\text{$\Tails\times\Zs$} \ar[r]^(0.55){\clutch_{(C/B,x)}}
&
\text{$\Tails$}
}
$$
\epr
The proof is similar to that of Proposition \ref{sing}, so we omit it.

\section{The intersection ring of $\Moz$}
\label{section ring of M02}

Just like $\Ws$, the Artin stack $\Moz=\BGm\times\Ws$ is not of finite type,
but it is exhaustible in the sense of \S\ref{section exhaustible}, with
exhaustion $(\BGm\times\Ws_n\subset\Moz)_n$. By Corollary \ref{A^*(XxW)}
we have
$$
A^*(\Moz)\isomorph R[t_0]
$$
where $t_0$ is transcendental over $R$.
We are interested in calculating the pull back morphism 
$\sigma^*:A^*(\Moz)\isomto A^*(\Moz)$
(in the sense of
\S\ref{section exhaustible})
induced by the involution
$\sigma:\Moz\isomto\Moz$
from the previous paragraph.

\bpr
\label{iota}
\barabic
\item
There is a unique ring homomorphism
$\iota:R[t_0]\to R[t_0]$ which 
maps the element $t_0$ to the element $-(t_0+t_1)$ and the element
$f(t_1,t_2,\dots,t_n)\in \Ic_n\subset R$ to the element $f(t_n,t_{n-1},\dots,t_1)$.
\item
The following diagram commutes
$$
\xymatrix@R=2ex{
\text{$A^*(\Moz)$} \ar[d]^{\isomorph}\ar[r]^{\sigma^*} 
&
\text{$A^*(\Moz)$} \ar[d]_{\isomorph} 
\\
R[t_0] \ar[r]^\iota
&
R[t_0]
}
$$ 
\earabic
\epr

\bpf
Let
$\iota_n:\Q[t_0,\dots,t_n]\to\Q[t_0,\dots,t_n]$
be defined by $\iota_n(t_0)=-(t_0+t_1+\dots+t_n)$ and
$\iota_n(t_i)=t_{n-i+1}$ for $1\leq i\leq n$. Observe that
$\iota_n$ maps the subring $R_n[t_0]\subset\Q[t_0,\dots,t_n]$ to itself.
Furthermore from the explicit description of $\sigma_n$ on page
\pageref{sigma} we see that the following diagram 
commutes:
$$
\xymatrix@R=2ex{
\text{$A^*(\BGm\times\Ws_n)$} \ar[d]^{\isomorph}\ar[r]^{\sigma_n^*} 
&
\text{$A^*(\BGm\times\Ws_n)$} \ar[d]_{\isomorph} 
\\
R_n[t_0] \ar[r]^{\iota_n}
&
R_n[t_0]
}
$$ 
The proposition now follows from the identity $\iota=\varprojlim_n\iota_n$,
which is easy to verify.
\epf

\brem
In a similar way we derive a formula for the pull-back morphism
$\clutch^*:A^*(\Moz)\to A^*(\Moz)\tensor A^*(\Moz)$.
However we omit it, since we will not need it in the sequel.

\erem
For the convenience of the reader let us review the structure of $R[t_0]$.
By definition we have
$$
R[t_0]=\Oplus_{n\geq 0}\left(\prod_{i=1}^nt_i\right)\Q[t_0,\dots,t_n]
\quad.
$$
Thus as a $\Q$-vector space, $R[t_0]$ has a basis consisting
of monomials
$$
t_0^{a_0}\prod_{i=1}^nt_i^{a_i}
\qquad
\text{where $n\geq 0$, $a_0\geq 0$, and $a_i\geq 1$ for $1\leq i\leq n$.}
$$
Let us denote by $x\tensor y\mapsto x\cdot y$ the multiplication 
on $R[t_0]$ which makes the natural morphism
$A^*(\Moz)\isomorph R[t_0]$ an isomorphism of rings.
On generators, this multiplication is given by
$$
\left(t_0^{a_0}\prod_{i=1}^nt_i^{a_i}\right)
\cdot
\left(t_0^{b_0}\prod_{i=1}^mt_i^{b_i}\right)
=
\sum_{k=0}^{m+n}
\sum_{\alpha,\beta}
t_0^{a_0+b_0}
\left(\prod_{i=1}^kt_i^{c_i}\right)
$$
where $(\alpha,\beta)$ runs through 
$\Hom_\Pc([1,n],[1,k])\times\Hom_\Pc([1,m],[1,k])$
with $\alpha([1,n])\cup\beta([1,m])=[1,k]$, and
$c_i\geq 1$ is defined by
$$
c_i=
\begin{cases}
a_j &\text{if $\alpha(j)=i$ for some $j\in[1,n]$ but $i\not\in\beta([1,m])$,}
\\
b_j &\text{if $\beta(j)=i$ for some $j\in[1,m]$ but $i\not\in\alpha([1,n])$,}
\\
a_p+b_q 
&\text{if $\alpha(p)=i$ for some $p\in[1,n]$, 
      and $\beta(q)=i$ for some $q\in[1,m]$.}
\end{cases}
$$
For example we have 
$$
t_1\cdot t_1=t_1^2+2t_1t_2
\qquad\text{and}\qquad
t_1\cdot(t_1t_2)=t_1^2t_2+t_1t_2^2+3t_1t_2t_3
\quad.
$$
The ring $(R[t_0],+,\cdot)$ is a commutative associative ring with unit element $1$.
The multiplication in this ring is compatible with the grading by the degree
of monomials: If $x,y\in R[t_0]$ are homogeneous of degree $d,e$, then
$x\cdot y$ is homogeneous of degree $d+e$.

On generators the involution $\iota:R[t_0]\isomto R[t_0]$ from
Proposition \ref{iota} is given by
$$
\iota\left(t_0^{a_0}\prod_{i=1}^nt_i^{a_i}\right)
=
(-t_0-t_1)^{\cdot a_0}\cdot\left(\prod_{i=1}^nt_i^{a_{n-i+1}}\right)
\quad.
$$
Notice that on the right hand side the power is taken with
respect to the multiplication ``$\cdot$''. For example we have
\begin{align*}
\iota(t_0^2t_1)=&(-t_0-t_1)\cdot(-t_0-t_1)\cdot t_1
\\
=&t_0^2t_1+2t_0(t_1\cdot t_1)+t_1\cdot t_1\cdot t_1
\\
=&t_0^2t_1+2t_0t_1^2+4t_0t_1t_2+t_1^3+2t_1^2t_2+2t_1t_2^2+6t_1t_2t_3
\quad.
\end{align*}
\section{The homology ring of $\Moz$}
\label{section homology of M02}

Proposition \ref{sing} implies in particular that the clutching
morphism $\clutch:\Zs\times\Zs\to\Zs$ is proper. 
If we take $(\Zs_n\subset\Zs)_n$ and $(\Zs_n\times\Zs_n\subset\Zs\times\Zs)_n$
as exhaustion for $\Zs$ and $\Zs\times\Zs$ respectively then $\clutch$
is compatible with these exhaustions. Thus by Lemma \ref{lem functorial}
there is a push forward morphism $\clutch_*:A_*(\Zs\times\Zs)\to A_*(\Zs)$.
From the isomorphism $\Zs\isomorph\BGm\times\Ws$ and Corollary \ref{A^*(XxW)}
it follows that $A_*(\Zs\times\Zs)=A_*(\Zs)\tensor A_*(\Zs)$.
Therefore $\clutch_*$ induces a bilinear map or multiplication
$$
A_*(\Zs)\tensor A_*(\Zs)\to A_*(\Zs)
$$          
which we denote by $x\tensor y\mapsto x\odot y$.
By Lemma \ref{Z monoid} this multiplication is associative
and $\sigma_*:A_*(\Zs)\isomto A_*(\Zs)$ is an anti-involution
of the ring $(A_*(\Zs),+,\odot)$.
The multiplication $\odot$ respects the grading of $A_*(\Zs)$:
For $x\in A_r(\Zs)$ and $y\in A_s(\Zs)$ we have $x\odot y\in A_{r+s}(\Zs)$.  
Since $A_0(\Zs)=(0)$ it follows that the ring $(A_*(\Zs),+,\odot)$
does not have a unity. We will see below that it is not commutative.

Since by Proposition \ref{M02} and \S\ref{section ring of M02}
we have $\Zs\isomorph\Moz$ and $A_*(\Moz)=A^{-1-*}\isomorph R[t_0]$,
the multiplication $\odot$ on $A_*(\Zs)$ induces a multiplication
on $A_*(\Moz)$ and on $R[t_0]$ which we denote by the same symbol.
Notice  that the isomorphism $A_*(\Moz)\isomorph R[t_0]$
involves a renumbering of the natural grading of $R[t_0]$:
A homogeneous polynomial 
$f\in \Ic_n[t_0]\subset R[t_0]$ 
of degree $r$ corresponds to an element in $A_{-1-r}(\Moz)$.
Thus for homogeneous elements $f,g\in R[t_0]$ of degree
$r,s$ respectively, the product $f\odot g$ is homogeneous of degree
$r+s+1$.

The following result allows us to explicitly calculate the $\odot$-product
of arbitrary elements in $R[t_0]$.

\bpr
\label{odot formula}
Let 
$$
f=t_0^{a_0}\prod_{i=1}^nt_i^{a_i}
\qquad\text{and}\qquad
g=t_0^{b_0}\prod_{i=1}^mt_i^{b_i}
$$
be two generators of $R[t_0]$.
Then for every $k\geq 1$
the image of $f\odot g$ by the canonical morphism $R[t_0]\to R_k[t_0]$
is
\begin{multline*}
Q_k(f,g):=
\\
\sum_{j=1}^k
\left(
\sum_{\alpha:[1,n]\injto[1,j-1]} f(t_0,t_{\alpha(1)},\dots,t_{\alpha(n)})
\right)
t_j
\left(
\sum_{\beta:[1,m]\injto[j+1,k]} 
g\left(\sum_{i=0}^jt_i,t_{\beta(1)},\dots,t_{\beta(m)}\right)
\right)
\end{multline*}
where $\alpha,\beta$ run through all strictly increasing maps 
as indicated below the sum signs.
\epr

An important ingredient for the proof of Proposition \ref{odot formula}
is the following lemma
which is an easy consequence of 
Lemma \ref{usefull}.

\blem
\label{ingredient}
Let $S$ be a scheme, let $L_0,\dots,L_k$ be line bundles on $S$ 
and for $1\leq i\leq k$ let $\lambda_i\in\Gamma(S,L_i)$
be global sections. Let $1\leq j\leq k$ and let 
$Y\subset S$ be the closed subscheme defined by the
vanishing of the section $\lambda_j$. Let $L'_i:=L_i|_Y$,
$\lambda'_i:=\lambda_i|_Y$.
Define
\begin{align*}
\xi&:=
|L_0,\Oo|\vdash(L_1,\lambda_1)\vdash\dots\vdash(L_k,\lambda_k)
&\in \Zs(S)
\quad,
\\
\xi'_1&:=
|L'_0,\Oo|\vdash(L'_1,\lambda'_1)\vdash\dots\vdash(L'_{j-1},\lambda'_{j-1})
&\in \Zs(Y)
\quad,
\\
\xi'_2&:=
|M',\Oo|\vdash(L'_{j+1},\lambda'_{j+1})\vdash\dots
        \vdash(L'_k,\lambda'_k)
&\in \Zs(Y)
\quad,
\end{align*}
where $M':=\Tensor_{i=0}^jL'_i$.
Then the restriction of the zollstock $\xi$ to the
closed subscheme $Y$ is canonically isomorphic to
the clutching of the two zollstocks $\xi'_1$ and $\xi'_2$.
I.e. we have
$$
\xi|_Y=\clutch(\xi'_1,\xi'_2)
\quad.
$$
\elem

\bpf
(of Proposition \ref{odot formula})
Let $\As:=\BGm\times[\A^1/\Gm]^k$, and
let $\xi$ be the \'etale morphism
$$
\As\to\BGm\times\Ws_k\isomto\Zs_k\injto\Zs
\quad.
$$ 
By definition, $\xi$ ``is'' the zollstock
$$
(X\to\As,x_1,x_2):=
|L_0,\Oo|\vdash(L_1,\lambda_1)\vdash\dots\vdash(L_k,\lambda_k)
$$
over $\As$, where 
$(L_0,(L_1,\lambda_1),\dots,(L_k,\lambda_k))$
is the universal object over $\As$. 
Let $\Ss$ be defined by the following Cartesian diagram:
$$
\xymatrix{
\text{$\Ss$}\ar[d]_{\clutch'}\ar[r]^(.45){\xi'} &
\text{$\Zs\times\Zs$} \ar[d]^\clutch
\\
\text{$\As$}\ar[r]^\xi &
\text{$\Zs$}
}
$$
By Proposition \ref{sing} we have $\Ss=\Sing(X/\As)$.
It follows directly from the definitions that $\Sing(X/\As)$
is nothing else but the disjoint union of the vanishing loci
of the universal sections $\lambda_j$. Therefore we have
$$
\Ss=\bigsqcup_{j=1}^kY_j
$$
where $Y_j$ is the closed substack
$$
\BGm\times[\A^1/\Gm]^{j-1}\times\BGm\times[\A^1/\Gm]^{k-j}
\quad\subset\quad\As
\quad.
$$
Let 
$$
(L^{(j)}_0,
(L^{(j)}_1,\lambda^{(j)}_1),
\dots,
(L^{(j)}_{j-1},\lambda^{(j)}_{j-1}),
L^{(j)}_j,
(L^{(j)}_{j+1},\lambda^{(j)}_{j+1}),
\dots,
(L^{(j)}_k,\lambda^{(j)}_k))
$$
be the universal object over $Y_j$, which is also the restriction
to $Y_j$ of the universal object over $\As$.
Now by Lemma \ref{ingredient} the morphism $Y_j\injto\Ss\to\Zs\times\Zs$
is given by a pair $(\xi^{(j)}_1,\xi^{(j)}_2)$ of zollstocks where
\begin{align*}
\xi^{(j)}_1&:=
|L^{(j)}_0,\Oo|\vdash(L^{(j)}_1,\lambda^{(j)}_1)\vdash\dots\vdash(L^{(j)}_{j-1},\lambda^{(j)}_{j-1})
&\in \Zs(Y_j)
\quad,
\\
\xi^{(j)}_2&:=
|M^{(j)},\Oo|\vdash(L^{(j)}_{j+1},\lambda^{(j)}_{j+1})\vdash\dots
        \vdash(L^{(j)}_k,\lambda^{(j)}_k)
&\in \Zs(Y_j)
\quad,
\end{align*}
and 
$M^{(j)}:=\Tensor_{i=0}^jL^{(j)}_i$.
Therefore identifying $A_*(Y_j)$ with $\Q[t_0,\dots,t_k]$, we have
\begin{multline*}
(\xi^{(j)}_1,\xi^{(j)}_2)^*(f,g)=
\\
\left(
\sum_{\alpha:[1,n]\injto[1,j-1]} f(t_0,t_{\alpha(1)},\dots,t_{\alpha(n)})
\right)
\left(
\sum_{\beta:[1,m]\injto[j+1,k]} 
g\left(\sum_{i=0}^jt_i,t_{\beta(1)},\dots,t_{\beta(m)}\right)
\right)
\end{multline*}
Since by \cite{Fulton} Theorem 3.3 and Example 3.3.2 (cf. also
the proof of Lemma \ref{lem2})
the push forward morphism 
$$
\Q[t_0,\dots,t_k]=A_*(Y_j)\to A_*(\As)=\Q[t_0,\dots,t_k]
$$
associated to the inclusion $Y_j\injto\As$
is nothing else but multiplication with $t_j$, it follows that
$\clutch'_*(\xi')^*(f,g)$ equals the expression $Q_k(f,g)$
stated in the Proposition. Now by Lemma \ref{lem functorial} (c)
we have $\clutch'_*(\xi')^*=\xi^*\clutch_*$.
Furthermore the pull-back
morphism
$$
R[t_0]=A_*(\BGm\times\Ws_k)\to A_*(\As)=\Q[t_0,\dots,t_k]
$$
is injective by Lemma \ref{lem5}. Therefore the Proposition follows.
\epf

As a first application of Proposition \ref{odot formula}
let us compute the product $f\odot g$
in the case where $g$ does not depend on $t_0$.
More precisely let $f\in\Ic_n[t_0]\subset R[t_0]$ 
and let $g\in\Ic_m\subset R$. We claim that in this case
we have
$$
f\odot g=f(t_0,\dots,t_n)t_{n+1}g(t_{n+2},\dots,t_{n+m+1})
\quad\in\quad\Ic_{n+m+1}[t_0]\subset R[t_0]
\quad.
$$
Indeed it suffices to prove this formula for monomials
$f=t_0^{a_0}\prod_{i=1}^nt_i^{a_i}$
and
$g=\prod_{i=1}^mt_i^{b_i}$.
Let $h$ be the element
$f(t_0,\dots,t_n)t_{n+1}g(t_{n+2},\dots,t_{n+m+1})$
of $\Ic_{n+m+1}[t_0]\subset R[t_0]$.
For any $k\geq 1$ the image of $h$ in $R_k[t_0]$ is
 (cf. proof of Theorem \ref{thm2})
$$
\sum_{\alpha:[1,n+m+1]\injto[1,k]}
h(t_0,t_{\alpha(1)},\dots,t_{\alpha(n+m+1)})
\quad,  
$$ 
where $\alpha$ runs through all strictly increasing maps
as indicated. But this expression is clearly equal to
the Polynomial $Q_k(f,g)$. Since $R[t_0]=\varprojlim'_kR_k[t_0]$,
the claim follows.
As a special case, let us compute the $n$-fold $\odot$-product
of the element $1\in R[t_0]$ with itself:
$$
1^{\odot n}=\prod_{i=1}^{n-1}t_i
\quad.
$$

It is more difficult to derive a formula for
$f\odot g$, if both $f$ and $g$ depend on $t_0$.
In this case it is often helpful to make use of the anti-involution $\iota$
as follows.
Let $f_1\in \Ic_n$, $g_1\in\Ic_m$.
Then we have 
$$
(t_0^af_1)\odot(t_0^bg_1)=
t_0^a(f_1\odot(t_0^bg_1))=
t_0^a\iota(\iota(t_0^bg_1)\odot\iota(f_1))=
t_0^a\iota(\iota(t_0^bg_1)\odot\iota(f_1))
\quad.
$$
Since the element 
$\iota(f_1)$ does not depend on $t_0$, calculation of the 
last $\odot$-product is easy. Notice however that the computation
of $\iota$ of an element that depends on $t_0$ can be quite 
involved (cf. example at the end of \S\ref{section ring of M02}).

As an application we will now prove a formula which we state as
a Lemma since it will be of use in \S\ref{section gravitation}.

\blem
\label{identity}
For any $n\geq 0$ we have the following identity in $R[t_0]$:
$$
1\odot(t_0+t_1)^{\cdot n}=t_1\cdot(t_0+t_1)^{\cdot n}
\quad.
$$
(The dot in the exponent means that we take the $n$-fold product
in the ring $(R[t_0],+,\cdot)$, not in the ring $\Q[t_0,t_1]$.)
\elem

\bpf
We exploit the fact that $\iota$ is an anti-involution
for $(R[t_0],+,\odot)$ and an involution for $(R[t_0],+,\cdot)$:
\begin{multline*}
1\odot(t_0+t_1)^{\cdot n}=
\iota((\iota(t_0+t_1))^{\cdot n}\odot 1)=
\iota((-t_0)^{n}\odot 1)=
\\
=
(-1)^n\iota(t_0^nt_1)=
(-1)^n(-t_0-t_1)^{\cdot n}\cdot t_1=
t_1\cdot(t_0+t_1)^{\cdot n}
\quad.
\end{multline*}
\epf


In Theorem \ref{structure} below we will show that the ring $(R[t_0],+,\odot)$
has a very simple structure. The following two lemmas will be useful
for its proof.

\blem
\label{dim}
The dimension of the subspace $R[t_0]_n\subset R[t_0]$ spanned by
monomials of degree $n$ equals $2^n$. 
\elem

\bpf
This follows by a simple calculation which we leave to the reader.
\epf

\blem
\label{1 odot x}
For all $x\in R[t_0]\setminus\{0\}$ we have $1\odot x\not\in t_0R[t_0]$.
\elem

\bpf
It suffices to show that for $x=t_0^r\prod_{i=1}^kt_i^{a_i}$ we have
$$
1\odot x=
t_1^{r+1}t_2^{a_1}\dots t_{k+1}^{a_k}
+ \text{(terms with strictly smaller power of $t_1$)}
\eqno{(*)}
$$
Indeed, from this equation it follows that for arbitrary nonzero $x$
the expansion of $1\odot x$ contains a monomial which is not divisible by $t_0$.
To prove $(*)$,
we make use of the identity $1\odot x=\iota(\iota(x)\odot 1)$.
We write
\begin{multline*}
\iota(x)\odot 1=
\left((-t_0-t_1)^{\cdot r}\cdot\prod_{i=1}^kt_{k-i+1}^{a_i}\right)\odot 1
=\\=
(-t_0)^r\left(\prod_{i=1}^kt_{k-i+1}^{a_i}\right)t_{k+1}
+ \text{(terms with strictly smaller power of $t_0$)}
\end{multline*}
Applying $\iota$ to the first term we get
\begin{multline*}
\iota\left(
(-t_0)^r\left(\prod_{i=1}^kt_{k-i+1}^{a_i}\right)t_{k+1}
\right)=
(t_0+t_1)^{\cdot r}\cdot\left(t_1\prod_{i=2}^{k+1}t_{i}^{a_{i-1}}\right)
=\\=
t_1^{r+1}\prod_{i=2}^{k+1}t_{i}^{a_{i-1}}
+ \text{(terms with strictly smaller power of $t_1$)}
\end{multline*}
Then it is clear that 
the first term in the last expression is the
term with the highest power of $t_1$ in the expansion of $1\odot x$.
This proves equation $(*)$.
\epf

\bthm
\label{structure}
As a $\Q$-vector space $R[t_0]$ has a basis consisting of 
elements of the form
$$
t_0^{a_1}\odot t_0^{a_2}\odot\dots\odot t_0^{a_r}
$$
where $r\geq 1$ and $a_i\geq 0$.
\ethm

\bpf
Observe first of all that the map $R[t_0]_n\to R[t_0]_{n+1}$,
$x\mapsto t_0x$ is injective. By Lemma \ref{1 odot x} the
map $R[t_0]_n\to R[t_0]_{n+1}$, $x\to 1\odot x$ is also injective
and the intersection of  $1\odot R[t_0]$ and $t_0\odot R[t_0]$
is trivial. By Lemma \ref{dim} it follows that
$$
R[t_0]_{n+1}=\left(1\odot R[t_0]_n\right)\oplus t_0R[t_0]_n
\quad.
$$ 
Observe that if $x\in R[t_0]_n$ is of the form 
$t_0^{a_1}\odot\dots\odot t_0^{a_r}$ then so is $1\odot x$ and $t_0x$.
Therefore the above equality shows that $R[t_0]_{n+1}$ has a basis
consisting of elements of the form 
$t_0^{a_1}\odot\dots\odot t_0^{a_r}$, if this is true for $R[t_0]_n$.
Obviously $R[t_0]_0$ has such a basis. So the theorem follows
by induction.
\epf

\bco
\label{basis}
As a $\Q$-vector space $R[t_0]$ has a basis consisting of 
elements of the form
$$
(t_0+t_1)^{a_1}\odot (t_0+t_1)^{a_2}\odot\dots\odot (t_0+t_1)^{a_r}
$$
where $r\geq 1$ and $a_i\geq 0$.
\eco

\bpf
By Proposition \ref{iota} the anti-involution $\iota$ maps $t_0$ to
$-t_0-t_1$. So the corollary follows directly from Theorem
\ref{structure}.
\epf

Let $B$ be a scheme (or stack) and let $(C/B,x)$ be a pointed curve over $B$.
Recall from \S\ref{section zollstock} p. \pageref{Z acts} that 
the monoid stack $Z$ acts on $\Tails:=\Tails(C/B,x)$ and that
by Proposition \ref{singt} the action morphism 
$\clutch_{(C/B,x)}:\Tails\times\Zs\to\Tails$
is proper. Therefore the push-forward with respect to this morphism
is defined and provides $A_*(\Tails)$ with the structure of 
a right module over the ring $(A_*(\Zs),+,\odot)$.

By the commutative diagram on page \pageref{cK diag} we see
that this module structure translates into a right $(R[t_0],+,\odot)$-module
structure on $A_*(B)\tensor R$, given by the push-forward with respect
to $\clutch_K:(B\times\Ws)\times\Moz\to B\times\Ws$
where $K$ is the line bundle $x^*\Oo(x)$ on $B$.

\bpr
\label{cK*}
We have
$$
(\clutch_K)_*:
\begin{cases}
\Bigl(A_*(B)\tensor R\Bigr)\tensor R[t_0]
& \to 
(A_*(B)\tensor R
\\
(a\tensor f)\tensor g
& \mapsto 
a(f\odot g)|_{t_0=-c_1(K)}
\end{cases}
$$
where for an element $h=\sum_ih_it_0^i\in R[t_0]$ with $h_i\in R$
and an element $c\in A_*(B)$ we  write
$$
h|_{t_0=c}:=\sum_ic^i\tensor h_i\quad\in\quad A_*(B)\tensor R
\quad.
$$
\epr

\bpf
The proof is entirely analogous to the one of Proposition \ref{odot formula}.
We omit the details.
\epf

For an application of our considerations in \S\ref{section gravitation} we 
need to understand the relationship between the $\phi$-class and the $\psi$-class.
We will show now that this relationship 
can conveniently be expressed as an identity in the ring $A_*(B)\tensor R$.
Let $(\Xs/\Tails,y)\to(C_\Tails/\Tails,x_\Tails)$ be the universal modification
over $\Tails$. Let $\psi\in A_*(\Tails)$ be the class of the
line bundle $y^*\omega_{\Xs/B}=y^*\Oo_{\Xs}(-y)$. By abuse of notation let $\psi$
denote also the image of that class by the canonical isomorphism
$A_*(\Tails)\isomto A_*(B)\tensor R$. On the other hand
let $\phi\in A_*(B)$ be the class of the line bundle 
$x^*\omega_{C/B}=x^*\Oo(-x)$ on $B$. 

\blem
\label{psi}
The following identity holds in $A_*(B)\tensor R$:
$$
\psi=\phi + t_1
\quad.
$$
\elem

\bpf
Let $(C_1/S,x_1)\to(C_2/S,x_2)$ be a simple modification.
By Proposition \ref{T1 isom A1/Gm} this modification corresponds to
a line bundle with section $(L,\lambda)$ on $S$.
More precisely it follows from the proof of this Proposition
that the line bundle $L$ is isomorphic to
$$
x_1^*\Oo_{C_1}(-x_1)\tensor x_2^*\Oo_{C_2}(x_2)
\quad.
$$
Now let $S$ be a $B$-scheme and let $(C',x')\to(C_S,x_S)$
be a modification over $S$ which via the isomorphism in
Prop. \ref{BxW isom T} corresponds to an object $(L_1,\lambda_1)\dots(L_r,\lambda_r)$
of $B\times\Ws$ over $S$. Since the pairs $(L_i,\lambda_i)$ 
come from decomposing $(C',x')\to(C_S,x_S)$ into a sequence of simple 
modifications and
a successive application of Prop. \ref{T1 isom A1/Gm}, it follows
that 
$$
\Tensor_iL_i=(x')^*\Oo_{C'}(-x')\tensor x_S^*\Oo_{C_S}(x_S)
\quad.
$$
From this the statement of the Lemma is immediate.
\epf

\bco
\label{cor psi}
Let $d\geq 1$ and $e\geq 0$.
The following formula holds in $A_*B\tensor R$:
$$
\psi^d\phi^e=\psi^{d-1}\phi^{e+1}+
\clutch_*(\phi^e\tensor(t_0+t_1)^{\cdot(d-1)})
\quad.
$$ 
where $\clutch:=\clutch_K:(B\times\Ws)\times\Moz\to B\times\Ws$ 
and $K=x^*\omega_{C/B}$.
\eco

\bpf
By Lemma \ref{psi} we have
$$
\psi^d\phi^e=\psi^{d-1}\phi^{e+1}+\phi^et_1(\phi+t_1)^{d-1}
$$
and by Proposition \ref{cK*} we have
$$
\clutch_*(\phi^e\tensor(t_0+t_1)^{\cdot(d-1)})=
\phi^e(1\odot(t_0+t_1)^{\cdot(d-1)})|_{t_0=\phi}
\quad.
$$
Therefore the Corollary follows from Lemma \ref{identity}.
\epf
\section{The extended operad}
\label{section operad}

Let $\Sb_n$ and $\Sb_{n+}$ denote the group of bijections of the
set $\{1,\dots,n\}$ and of the set $\{0,\dots,n\}$ respectively.
We consider $\Sb_n$ as a subgroup of $\Sb_{n+}$ in the obvious way.
It is well know that
the group $\Sb_{n+}$ is generated by the elements of $\Sb_n$ together
with the
$(n+1)$-cycle  $\tau_n:=(0\dots n)\in\Sb_{n+}$.

\bdefin
(Cf. \cite{Getzler} \S1.5)
A {\em cyclic operad} $\Pc$ (without unity, 
in the category of $\Q$-vector spaces) 
is a sequence $(\Pc(n))_{n\geq 1}$
of $\Q$-vector spaces together with  the following structures:
\balph
\item
For every $n\geq 1$ there is given a left action of the group
$\Sb_{n+}$ on $\Pc(n)$.
\item
For every $m,n\geq 1$ and $1\leq j\leq m$ there is given 
a composition morphism
$$
\comp_j:\Pc(m)\tensor\Pc(n)\to\Pc(m+n-1)
$$
\ealph
These structures are required to satisfy the following axioms
\barabic
\item
For $\pi\in\Sb_m$, $\rho\in\Sb_n$  and $1\leq j\leq m$ let
$\pi\comp_j\rho\in\Sb_{m+n-1}$ be defined as follows:
Set
\begin{align*}
(u_1,\dots,u_{j-1},u_{j+1},\dots,u_m)&:=(1,\dots,j-1,j+n,\dots,m+n-1)
\quad,
\\
(u'_1,\dots,u'_n)&:=(j,j+1,\dots,j+n-1)
\quad.
\end{align*}
Then the inverse of $\pi\comp_j\rho$ is defined by
$$
(\pi\comp_j\rho)^\mo(i):=
\begin{cases}
u_{\pi^\mo(i)} &
\quad\text{for $1\leq i\leq \pi(j)-1$}
\\
u'_{\rho^\mo(i-\pi(j)+1)} &
\quad\text{for $\pi(j)\leq i\leq \pi(j)+n-1$}
\\
u_{\pi^\mo(i-n+1)} &
\quad\text{for $\pi(j)+n\leq i\leq m+n-1$}
\end{cases}
$$
With this notation, for any $m,n\geq 1$, $1\leq j\leq m$,
$\pi\in\Sb_m$, $\rho\in\Sb_n$, 
$a\in\Pc(m)$ and $b\in\Pc(n)$
the following equality holds
$$
(\pi\comp_j\rho)(a\comp_j b)=(\pi a)\comp_{\pi(j)}(\rho b)
\quad.
$$
\item
For any $m,n\geq 1$ and for any
$a\in\Pc(m)$ and $b\in\Pc(n)$
the following equality holds
$$
\tau_{m+n-1}(a\comp_m b)=(\tau_nb)\comp_1(\tau_ma)
\quad.
$$
\item
For any $k,l,m\geq 1$ and $1\leq i<j\leq k$ 
and for any $a\in\Pc(k)$, $b\in\Pc(l)$, $c\in\Pc(m)$:
$$
(a\comp_ib)\comp_{j+l-1}c=(a\comp_jc)\comp_ib
\quad.
$$
\item
For any $k,l,m\geq 1$, $1\leq i\leq k$, $1\leq j\leq l$
and for any $a\in\Pc(k)$, $b\in\Pc(l)$, $c\in\Pc(m)$:
$$
(a\comp_ib)\comp_{i+j-1}c=a\comp_i(b\comp_jc)
\quad.
$$
\earabic
\edefin

\brem
\barabic
\item
In \cite{Getzler} (1.2.1) there is an error in 
the explicit definition of the permutation
$\pi\comp_j\rho$.  
The axioms can be visualised by drawing trees.
We refer to \cite{Getzler} \S1 for further motivation.
\item
It is easy to express axioms (1)-(4) in terms of commutative
diagrams. By means of these diagrams one can define cyclic operads
in any monoidal category.
\earabic
\erem

\bex
\label{EMb}
Let $\kb$ be a ring.
Let $M$ be an $\kb$-module and let $b:M\tensor M\to \kb$
be a bilinear form.
We define a cyclic operad $\E[M,b]$ by putting
$
\E[M,b](n):=M^{\tensor(n+1)}
$,
providing this space with the natural left $\Sb_{n+}$-action 
$$
\pi(v_0\tensor\dots\tensor v_n)
:=v_{\pi^\mo(0)}\tensor\dots\tensor v_{\pi^\mo(n)}
\quad,
$$
and 
with the composition morphisms
$$
\comp_j:
M^{\tensor(m+1)}\tensor M^{\tensor(n+1)}\to M^{\tensor(m+n)}
$$
which map an element
$(v_0\tensor\dots\tensor v_m)\tensor(w_0\tensor\dots\tensor w_n)$
to the element
$$
b(v_j\tensor w_0)
v_0\tensor\dots\tensor v_{j-1}\tensor
w_1\tensor\dots\tensor w_n\tensor
v_{j+1}\tensor\dots\tensor v_m
\quad.
$$
\eex

\bex
\label{EsMb}
There is an obvious super version of Example \ref{EMb} which
will play an important role in \S\ref{section gravitation}:
Let $\kb$ be a ring and let $M=M_0\oplus M_1$ be a $\Z/2$-graded 
$\kb$-module. For a homogeneous element $v\in M$ we denote
by $|v|\in\{0,1\}$ its degree.
Let $b:M\tensor M\to\kb$ be a symmetric 
bilinear
form of degree $0$.
This means that for homogeneous elements $v,w\in M$
we have $b(v,w)=(-1)^{|v||w|}b(w,v)$ and $b(v,w)=0$ if
$|v|\neq|w|$. 
We define the cyclic operad $\E^s[M,b]$ as follows.
The objects of the operad are
$
\E^s[M,b](n):=M^{\tensor(n+1)}
$.
The $\Sb_{n+}$-action on $M^{\tensor(n+1)}$ is given by
$$
\pi(v_0\tensor\dots\tensor v_n)
:=
(-1)^{|v_j||v_{j+1}|}
v_0\tensor\dots\tensor v_{j-1}\tensor v_{j+1}\tensor v_j
\tensor v_{j+2}\tensor\dots\tensor v_n
$$
where the $v_i$ are homogeneous elements of
$M$ and $\pi=(j,j+1)$ is the transposition of  two neighbouring
elements $j$ and $j+1$.
Finally, for homogeneous elements $v_i$, $w_i$ the composition
$
(v_0\tensor\dots\tensor v_m)\comp_j(w_0\tensor\dots\tensor w_n)
$
is defined as
$$
(-1)^N
b(v_j, w_0)
v_0\tensor\dots\tensor v_{j-1}\tensor
w_1\tensor\dots\tensor w_n\tensor
v_{j+1}\tensor\dots\tensor v_m
$$
where $N=\sum_{r=0}^n\sum_{s=j+1}^m|w_r||v_s|$
\eex

\bex
\label{Mb0}
The cyclic operad $\Mb_0$ in the monoidal 
category of algebraic varieties
is defined as follows:
$\Mb_0(1)$ is the empty variety and for $n\geq 2$ the variety
$\Mb_0(n)$  is the moduli space $\Mb_{0,n+1}$ of stable $(n+1)$-pointed
curves of genus zero. This space is endowed with an obvious left
$\Sb_{n+}$-action. For $m,n\geq 2$ the composition morphism 
$$
\comp_j:\Mb_{0,m+1}\times\Mb_{0,n+1}\to\Mb_{0,n+m}
$$
is the clutching morphism, which in terms of $S$-valued points
can be expressed by
$$
\Bigl(\bigl(C,(x_i)\bigr),
\bigl(C',(x'_i)\bigr)\Bigr)
 \mapsto 
\Bigl(
C\sqcup C'/(x_j=x'_0); x_0,\dots,x_{j-1},x'_1,\dots,x'_n,x_{j+1},\dots,x_m
\Bigr)
$$
\eex

\bex
\label{AMb0}
By applying the homology functor $A_*$ to the operad $\Mb_0$, we
obtain a cyclic operad $A\Mb_0$ in the category of $\Q$-vector spaces.
\eex

\bex
\label{PxQ}
Let $\Pc$ and $\Qc$ be two cyclic
operads. Then the product cyclic operad $\Pc\tensor\Qc$
has objects $(\Pc\tensor\Qc)(n):=\Pc(n)\tensor\Qc(n)$,
with $\Sb_{n+}$ acting on both factors, and composition
$$
(a\tensor a')\comp_j(b\tensor b'):=(a\comp_j b)\tensor(a'\comp_j b')
$$
for $a\in\Pc(m)$, $a'\in\Qc(m)$, $b\in\Pc(n)$, $b'\in\Qc(n)$. 
\eex

\vspace{2mm}
The central object of this paper is the cyclic operad $\Mt_0$ in the monoidal
category of Artin stacks, which we will define now. We let
\begin{align*}
\Mt_0(1):=&\Moz=\BGm\times\Ws \\
\Mt_0(n):=&\Mb_{0,n+1}\times\Ws^{n+1} \text{\quad for $n\geq 2$.}
\end{align*}
Thus for $n\geq 2$ an $S$-valued point of $\Mt_0(n)$ is a tupel
$$
(C,(x_0,\dots,x_n),(\xi_0,\dots,\xi_n))
$$
where $(C,(x_i))$ is an $n$-pointed stable curve of genus $0$ over $S$
and $\xi_i$ is an object in $\Ws(S)$ for $0\leq i\leq n$.

On $\Mt_0(1)$ the action of $\Sb_{1+}$ is given by the
involution $\sigma$ defined on page \pageref{clutch sigma}.
For $n\geq 2$ there is an obvious left action of the symmetric group $\Sb_{n+}$
on $\Mt_0(n)$ which on $S$-valued points is given by
$$
\tau(C,(x_0,\dots,x_n),(\xi_0,\dots,\xi_n)):=
(C,(x_{\tau(0)},\dots,x_{\tau(n)}),
(\xi_{\tau(0)}\dots,\xi_{\tau(n)}))
\quad\text{for $\tau\in \Sb_{n+}$.}
$$

Let us now define the composition morphisms
$$
\comp_j:
\Mt_0(m)\times\Mt_0(n)\to\Mt_0(m+n-1)
\quad\text{for $m,n\geq 1$ and $1\leq j\leq m$.}
$$
Assume first that
$m,n\geq 2$. Then we let $\comp_j$
be the morphism which on $S$-valued points is given by
$$
\bigl(C,(x_i),(\xi_i)\bigr)
\comp_j
\bigl(C',(x'_i),(\xi'_i)\bigr)
=
\bigl(C'',(x''_i),
(\xi_0,\dots,\xi_{j-1},\xi'_1,\dots,\xi'_n,\xi_{j+1},\dots,\xi_m)\bigr)
$$
where
$$
\bigl(C'',(x''_i)\bigr):=
\bigl(C,(x_i)\bigr)
\comp_j
\bigl(C',(x'_i)\bigr)
\quad,
$$
as defined in Example \ref{Mb0}.

Next assume that $m\geq 2$ and $n=1$. In this case we define the morphism
$
\comp_j:
\Mt_0(m)\times\Mt_0(1)\to\Mt_0(m)
$
by
$$
\bigl(C,(x_i),(\xi_i)\bigr)
\comp_j
\eta
:=
\bigl(C,(x_0,\dots,x_m),(\xi_0,\dots,\xi_{j-1},
\clutch_K(\xi_j,\eta),\xi_{j+1},\dots,\xi_m)\bigr)
$$
where $K:=x_j^*\Oo_{C}(x_j)$ and $\clutch_{K}$  is the morphism
from Definition \ref{cK def}.

Now let $m=1$ and $n\geq 2$. Then the morphism
$
\comp_1:
\Mt_0(1)\times\Mt_0(n)\to\Mt_0(n)
$
is defined by
$$
\eta
\comp_1
\bigl(C,(x_i),(\xi_i)\bigr)
:=
\bigl(C,(x_0,\dots,x_n),
(\clutch_K(\xi_0,\sigma(\eta)),\xi_1,\dots,\xi_n)\bigr)
$$
where $K:=x_0^*\Oo_{C}(x_0)$,
$\clutch_{K}$  is the morphism
from Definition \ref{cK def}, and $\sigma:\Mt_0(1)\isomto\Mt_0(1)$
is the involution from page \pageref{clutch sigma}.

Finally let $m=n=1$. Then we let
$
\comp_1:\Mt_0(1)\times\Mt_0(1)\to\Mt_0(1)
$
be the morphism $\clutch$ from page \pageref{clutch sigma}.

\vspace{1mm}
It is straightforward to check that $\Mt_0$ satisfies the
axioms of a cyclic operad (without unity) in the monoidal category
of Artin
stacks.
Since $\Mt_0(1)$ is non-empty and for $n\geq 2$ the objects
$\Mt_0(n)$ and $\Mb_0(n)$ (cf. Example \ref{Mb0}) 
differ only by a factor $\Ws^{n+1}$, 
it makes sense to call the operad $\Mt_0$ an {\em extension} of the operad
$\Mb_0$.

\vspace{3mm}
We would like to define a cyclic operad $A\Mt_0$
in the monoidal category of vector spaces
by applying the the cohomology functor $A_*$ to $\Mt_0$. 
However there is a problem here, since the morphisms 
$
\comp_j:
\Mt_0(m)\times\Mt_0(n)\to\Mt_0(m+n-1)
$
are not proper except in those cases where $m=1$ or $n=1$,
and we do not know how functoriality of homology can be 
extended to the morphisms $\comp_j$ also in the cases
where both $m,n\geq 2$.
Therefore in the following definition of the cyclic operad
$A\Mt_0$ the composition morphisms 
$A\Mt_0(m)\tensor A\Mt_0(n)\to A\Mt_0(m+n-1)$
for $m,n\geq 2$ will be defined in an {\em ad hoc} manner.

For $n\geq 1$ let $A\Mt_0(n):=A_*(\Mt_0(n))$.
By Corollary \ref{A^*(XxW)} we can identify these groups as follows:
\begin{align*}
A\Mt_0(1)=&R[t_0]\quad,\\
A\Mt_0(n)=& A_*\Mb_{0,n+1}\tensor R^{\tensor(n+1)}
\quad\text{for $n\geq 2$.}
\end{align*}

In the range $n\geq 2$ we define $A\Mt_0$ to have the same structure as
the product operad $A\Mb_0\tensor\E[R,b]$
(cf. Examples \ref{AMb0}, \ref{EMb}, \ref{PxQ})
where the bilinear form 
$b:R\tensor R\to\Q$ is defined as follows:
Let $\pr_0:R=\Oplus_{i=0}\Ic_i\to\Ic_0=\Q$ be the projection onto
the component of degree zero. Then we set
$$
b(x\tensor y):=\pr_0(x)\pr_0(y)
\quad\text{for $x,y\in R$.}
$$

It remains to define the structure of $A\Mt_0$ in the cases where
the first component $A\Mt_0(1)$ is involved. 
In these cases we can take the structure induced by the
operad $\Mt_0$. 
By the results of \S\ref{section homology of M02}
this can be made explicit as follows.

The action of $\Sb_{1+}$ on $R[t_0]$ is given by the involution
$\iota$ from Proposition \ref{iota}.

The composition 
morphism $\comp_1:A\Mt_0(1)\tensor A\Mt_0(1)\to A\Mt_0(1)$
is identified with the multiplication $\odot:R[t_0]\tensor R[t_0]\to R[t_0]$
(cf. Proposition \ref{odot formula}).

For $m\geq 2$ and $1\leq j\leq m$ the composition morphism 
$\comp_j:A\Mt_0(m)\tensor A\Mt_0(1)\to A\Mt_0(m)$
is identified with the morphism
$
(A_*\Mb_{0,m+1}\tensor R^{\tensor(m+1)})\tensor R[t_0]
\to 
A_*\Mb_{0,m+1}\tensor R^{\tensor(m+1)}
$
which maps an element
$(a\tensor (f_0\tensor\dots\tensor f_m))\tensor g$
to the element
$$
(a\tensor f_0\tensor\dots\tensor f_{j-1}\tensor 1\tensor
f_{j+1}\tensor\dots\tensor f_m)
q'_j\Bigl((f_j\odot g)|_{t_0=\phi_j}\Bigr)
$$
(cf. Proposition \ref{cK*}).
In the last formula, the symbol $\phi_j$ denotes the class
$$
c_1(x_j^*\omega_{m+1})
\quad\in\quad A_*\Mb_{0,m+1}
$$
where 
$(\Xs_{0,m+1}/\Mb_{0,m+1},(x_i))$ is the universal object over
$\Mb_{0,m+1}$, $\omega_{m+1}$ is the relative dualizing sheaf of
$\Xs_{0,m+1}\to\Mb_{0,m+1}$,
and $q'_j$ is the ring homomorphism
\label{q'}
$$
\begin{cases}
A_*\Mb_{0,m+1}\tensor R &\to A_*\Mb_{0,m+1}\tensor R^{\tensor(m+1)}\\
x\tensor f &\mapsto 
x\tensor 
\underset{\text{$j$ times}}
{\underbrace{1\tensor\dots\tensor 1}}
\tensor f\tensor 
\underset{\text{$(m-j)$ times}}
{\underbrace{1\tensor\dots\tensor 1}}
\end{cases}
$$

Finally, for $n\geq 2$ the composition morphism
$\comp_1:A\Mt_0(1)\tensor A\Mt_0(n)\to A\Mt_0(n)$
is identified with the morphism
$
R[t_0]\tensor (A_*\Mb_{0,n+1}\tensor R^{\tensor(n+1)})
\to 
A_*\Mb_{0,n+1}\tensor R^{\tensor(n+1)}
$
which maps an element
$g\tensor(a\tensor (f_0\tensor\dots\tensor f_n))$
to the element
$$
(a\tensor 1\tensor f_1\tensor\dots\tensor f_n)
q'_0\Bigl((f_0\odot \iota(g))|_{t_0=\phi_0}\Bigr)
\quad.
$$

\bpr
Let $\psit_1$ be the element $-t_0-t_1$ of $R[t_0]=A\Mt_0(1)$.
Let $n\geq 1$ and let $x\in A\Mt_0(n)$. There exists a 
sequence of elements
$x_0,x_1,\dots,x_r\in A\Mt_0(n)$, a sequence of numbers
$a_1,\dots,a_r\in\N_0$, and a sequence of numbers 
$j_1,\dots,j_r\in\{1,\dots,n\}$ such that
\broman
\item
$x_r$ equals $x$,
\item
for $1\leq i\leq r$ we have $x_i=x_{i-1}\comp_{j_i}\psit_1^{a_i}$,
\item
$x_0$ equals $\psit_1^{a_0}$ for some $a_0\geq 0$, if $n=1$,
and $x_0$ is an element of $A_*\Mb_{0,n+1}$, if $n\geq 2$.
\eroman
\epr

\bdefin
\label{morph}
Let $\Pc$ and $\Qc$ be two cyclic operads.
A {\em morphism} $$F:\Pc\to\Qc$$ is a sequence $F(n):\Pc(n)\to\Qc(n)$
($n\geq 1$) of morphisms  such that each $F(n)$ commutes with
the $\Sb_{n+}$-operation and such that for every $m,n\geq 1$, $1\leq j\leq m$,
$x\in\Pc(m)$, $y\in\Pc(n)$ we have
$$
F(m+n-1)(x\comp_j y)=F(m)(x)\comp_jF(n)(y)
\quad.
$$
\edefin

\bex
\label{AMb to AMt}
There is a natural morphism $F$ from the cyclic operad $A\Mb_0$
(cf. Example \ref{Mb0}) to the cyclic operad $A\Mt_0$ defined above.
$F(1)$ is the zero map and for $n\geq 2$ the morphism
$F(n):A\Mb_0(n)\to A\Mt_0(n)$ is defined
by
$$
F(n)(x):=x\tensor \underset{\text{$(n+1)$ times}}
{\underbrace{1\tensor\dots\tensor 1}}
\quad
\text{for $x\in\Mb_{0,n+1}$.}
$$ 
\eex

\section{The morphisms $\stt:\Mb_{g,n}(V,\beta)\to\Mt_{g,n}$}
\label{section stt}
Let $V$ be a smooth projective variety and let 
$\Num_1(V)\subset Z_1(V)$ be the subgroup of $1$-dimensional cycles which are
numerically equivalent to zero. It is well known (cf. \cite{Fulton}, Example
19.1.4) that the group
$N_1(V):=Z_1(V)/\Num_1(V)$ is a finitely generated free abelian group.
Let
$$
B(V):=\{\beta\in N_1(V) \st\text{$c_1(L)\cap\beta\geq 0$ for all ample line
bundles $L$ on $V$}\}
\quad.
$$
Then $B(V)$ is a submonoid of $N_1(V)$ with the property that
for all $\beta\in B(V)$ there are only finitely many pairs
$(\beta_1,\beta_2)\in B(V)\times B(V)$ such that $\beta_1+\beta_2=\beta$.

Choose a basis $(e_1,\dots,e_r)$ of $N_1(V)$,  and let
\label{Lambda} $\Lambda$ be the localisation by $\prod_{i=1}^rq_i$
of the power series ring $\Q[[q_1,\dots,q_r]]$.
Then the mapping
$$
\begin{cases}
B(V) &\to \Lambda \\
\beta &\mapsto q^\beta:=q_1^{b_1}\dots q_r^{b_r}
\quad\text{where $\beta=\sum_ib_ie_i$}
\end{cases}
$$
is a  character of the monoid $B(V)$.

For $g,n\geq 0$ and $\beta\in B(V)$ let 
$\Mb_{g,n}(V,\beta)$ be the moduli stack of stable $n$-pointed maps
of genus $g$ and class $\beta$ into $V$. We refer to \cite{Manin} Ch.V,\S5
for a precise definition of $\Mb_{g,n}(V,\beta)$.
Recall from \cite{Manin} Ch.V \S4.4 that for $(g,n)$ in the stable range
(i.e. for $n>2(1-g)$)
there exists a natural morphism
$$
\sta:\Mb_{g,n}(V,\beta)\to\Mb_{g,n}
$$
called the ``stabilisation morphism''.
For $(g,n)=(0,2)$ there is no such morphism, since $\Mb_{0,2}$
is the empty variety.

Let $\Moz:=\BGm\times\Ws$ as in Definition \ref{M02}, and for $(g,n)$ 
in the stable range let 
$$
\Mt_{g,n}:=\Mb_{g,n}\times \Ws^n
\quad.
$$
Thus the stacks $\Mt_{g,n}$ are well defined and non-empty for
all $(g,n)$ with $n>1-g$.
As a manner of speaking, we say that $g,n$ in the 
{\em half-stable range},
if the condition $n>1-g$ is fulfilled.
The following result shows that in the half-stable range
there is a morphism analogous to the morphism $\sta$.

\bpr
\label{stt}
Let $g,n$ be in the half-stable range. 
Then there exist natural morphisms
$$
\stt:\Mb_{g,n}(V,\beta)\to\Mt_{g,n}     
$$
such that in the stable range we have 
$\sta=\pr_1\comp\stt$ where $\pr_1:\Mt_{g,n}\to\Mb_{g,n}$ is the
projection onto the first factor.
\epr

For the proof we need some preparation.

\bdefin
Let $(g,n)$ be in the half-stable range.
Let $(C,(x_i))$ be a prestable $n$-pointed curve of genus $g$ 
over an algebraically closed field and let $x\in\{x_0,\dots,x_{n-1}\}$
be one of the marked points. 
\balph
\item
A {\em end} of $(C,(x_i))$ is a
closed subscheme of $C$ which is isomorphic to a chain of projective
lines, and
meets the rest of $C$ in exactly one point which lies 
on one of the extremal components (called the {\em base component}) 
of the chain.
\item
An end of of $(C,(x_i))$ is called a {\em loose end}  if it 
contains none of the marked points.
\item
An end of of $(C,(x_i))$ is called a {\em pointed end at $x$},  
if $x$ is contained in the end, if there is no other
marked point besides $x$ contained in the end, and if either the end 
consists of only one component,
or $x$ lies in the extremal component opposite to the base
component of the end.
\item
A morphism $(C,(x_i))\to(C',(x'_i))$ of prestable pointed curves
is called a 
{\em contraction of loose ends} (resp. {\em contraction of a pointed end at $x$}),
if for every point $y'\in C'$
the morphism $C\to C'$ is either an isomorphism at $y'$, or
its fibre over $y'$ is a loose end 
(resp. a pointed end at $x$) 
of $(C,(x_i))$.
Observe that in the pointed case the morphism
$f:(C,(x_i))\to(C',(x'_i))$ is what in \S\ref{section modifications}
we called a {\em modification} of $C'$ at $f(x)$.
\ealph
\edefin

\blem
\label{lem loose}
Let $(g,n)$ be in the half-stable range and 
let $(C/S,(x_i))$ be a prestable $n$-pointed curve of genus $g$ 
over a scheme $S$.
Assume that there is a line bundle  $L$ over $C$ such that 
\barabic
\item
the restriction of $L$ to any component of a geometric
fibre of $C\to S$ is of non-negative degree,
\item
the degree of the restriction of $L$ to a geometric
fibre of $C\to S$ does not depend on the fibre,
\item
the restriction of the sheaf 
$\omega_{C/S}(\sum_ix_i)\tensor L$ to any geometric
fibre of $C\to S$ is ample.
\earabic
Then there is a prestable pointed curve $(C'/S,(x'_i))$ 
whose geometric fibres are without loose ends, and a morphism
$$
f:(C/S,(x_i))\to(C'/S,(x'_i))
$$
which is fibre-wise a contraction of loose ends.
Furthermore there is a line bundle $L'$ on $C'$ which
satisfies properties (1)-(3) above.

If $f':(C/S,(x_i))\to(C''/S,(x''_i))$ is a second morphism with
these properties then there exists a unique isomorphism
$h:(C'/S,(x'_i))\isomto(C''/S,(x''_i))$
such that $f'=h\comp f$.
\elem

\bpf
Assume first that $S$ is the spectrum of an algebraically
closed field. Observe that for every irreducible component $Y$
of $C$ the degree of the restriction of $\omega_{C/S}(\sum_ix_i)$ to
$Y$ is either non-negative, or it is $-1$. The latter occurs
precisely if $Y$ is a loose end.
Therefore, if $(C/S,(x_i))$ has loose ends, then
the set of all $e\geq 1$ 
such that the restriction of $M_e:=\omega_{C/S}(\sum_ix_i)^e\tensor L$
to every component of $C$ is of non-negative degree, is bounded from above
by $\deg(L)$. Let $e$ be maximal with this property.
Then the restriction of $M_e$ to a component $Y$ of $C$ is of degree
zero if and only if $Y$ is a loose end.

There is an $m\geq 1$ which depends only on the numbers
$g$, $n$, and $\deg(L)$,
such that $H^1(C,M_e^m)$ vanishes, and $H^0(C,M_e^m)\tensor\Oo_C\to M_e^m$
is surjective. The latter defines a morphism
$$
C\to\Pp(H^0(C,M_e^m))
$$ 
and we let $(C^{(1)},(x^{(1)}_i))$ be the image of $(C,(x_i))$
by this morphism. Then $(C,(x_i))\to(C^{(1)},(x^{(1)}_i))$ is
a non-trivial contraction of loose ends.

Observe that $(C^{(1)},(x^{(1)}_i))$ is an $n$-pointed prestable
curve of genus $g$ and that the restriction $L^{(1)}$
of $\Oo_{\Pp}(1)$ to $C^{(1)}$ satisfies
properties (1) and (2) of the Lemma.
Therefore if $(C^{(1)},(x^{(1)}_i))$ has loose ends, we  can 
apply the same procedure
to this curve. 
Iterating this process we obtain a sequence
$$
(C,(x_i))\to(C^{(1)},(x^{(1)}_i))\to(C^{(2)},(x^{(2)}_i))\to\dots
$$
of non-trivial contractions of loose ends.
The length of this sequence is bounded by the number of components
of $C$, which in turn is bounded by a number which depends
only on $(g,n,\deg(L))$. Thus there is a last pointed curve 
$(C',(x'_i))$ in this sequence which is without loose ends,
and it is clear that the composed morphism $(C,(x_i))\to(C',(x'_i))$
is a contraction of loose ends.

Let $\Ll$ be the pull-back to $C$ of a very ample line bundle on $C'$.
Then we have 
$$
C'=\Proj\left(\oplus_{k=0}^\infty H^0(C,\Ll^k)\right)
$$
Furthermore, the line bundle
$\Ll$ has the property that its restriction to any component
is non-negative and that it is trivial if and only if the component
belongs to a loose end. For any other line bundle $\Ll'$ on
$C$ with this property the scheme
$\Proj\left(\oplus_{k=0}^\infty H^0(C,(\Ll')^k)\right)$
is isomorphic to $C'$. This proves that $C'$ is determined up to
isomorphism.

The restriction of the morphism $C\to C'$ to non-contracted
components is an isomorphism. Therefore there is no non-trivial
automorphism of $C'$ as an object under $C$.
This shows that $C'$ is determined even up to a {\em unique} isomorphism.

This finishes the proof in the case where $S$ is the spectrum
of an algebraically closed field.
The general case works essentially the same way, one just
has to replace
the cohomology groups by the appropriate direct image sheaves.
\epf

\blem
\label{lem contract}
Let $(g,n)$ be in the stable range, 
let $(C/S,(x_i))$ be a prestable $n$-pointed curve of genus $g$,
whose geometric fibres are without loose
ends, and let $x$ be one of the $x_i$.
Assume that there is a line bundle  $L$ over $C$ such that
the conditions (1)-(3) of Lemma \ref{lem loose} are satisfied. 
Then there is a morphism
$$
f:(C/S,(x_i))\to(C'/S,(x'_i))
$$
of prestable pointed curves,
such that the geometric fibres of $(C'/S,(x'_i))$ are
without a pointed end at $f(x)$,
and such that $(C/S,x)\to(C'/S,f(x))$ is a modification of $C'$ at $f(x)$. 
Furthermore there is a line bundle $L'$ on $C'$ which again
satisfies properties (1)-(3) stated in \ref{lem loose}.

If $f':(C/S,(x_i))\to(C''/S,(x''_i))$ is a second morphism with
these properties then there exists a unique isomorphism
$h:(C'/S,(x'_i))\isomto(C''/S,(x''_i))$
such that $f'=h\comp f$.
\elem

\bpf
As in the proof of Lemma \ref{lem loose} we first assume that
$S$ is the spectrum of an algebraically closed field.
Arguing similarly as in Lemma \ref{lem loose}, with the divisor 
$\sum_ix_i$ replaced by $\sum_ix_i-x$, we see that there is
a line bundle on $C$ whose restriction to every component
$Y$ is of non-negative degree, and is of degree zero if and only
if $Y$ is a pointed end at $x$.
Now the sought for morphism 
$$
(C,(x_i))\to (C',(x'_i))
$$
is constructed by an iterative process exactly as in
the proof of the previous lemma.

Uniqueness and the case of a general base scheme follows
also as in the proof of Lemma \ref{lem loose}.
\epf

After these preparations we are now able to prove Proposition
\ref{stt}.

\bpf(of Proposition \ref{stt})
Let $(C/S, (x_i), f:C\to V)$ be a stable $n$-pointed map of genus
$g$ and of class $\beta\in B(V)$ to $V$ over a scheme $S$. 

Let $M$ be a very ample line bundle on $V$.
It is well known (cf. \cite{Manin}, Ch.V, proof of Prop. 4.4.1)
that the sheaf $f^*(M)^{\tensor 3}$ satisfies
properties (1)-(3) stated in Lemma \ref{lem loose}.

Applying this Lemma, we obtain an (up to a unique isomorphism) unique
morphism
$$
(C/S,(x_i))\to(C'/S,(x'_i))
$$
such that the fibres of $(C'/S,(x'_i))$ are without
loose ends and a line bundle $L'$ on $C'$ which again
satisfies
properties (1)-(3) in Lemma \ref{lem loose}.

If $(g,n)=(0,2)$ then $(C'/S,x'_1,x'_2)$ is a zollstock over $S$
which by Proposition \ref{M02} gives rise to an object $\eta\in\Moz(S)$,
and we set $\stt(C/S, (x_i), f):=\eta$.

Assume now that $(g,n)$ is in the stable range.
Let $1\leq j\leq n$. By Lemma \ref{lem contract} there is an
(up to a unique isomorphism) unique modification $(C'/S,x'_j)\to(C''/S,x''_j)$
such that $C''$ is not further contractible at $x''_j$.
Let $\xi_j\in\Ws(S)$ be the object which by Proposition \ref{BxW isom T}
is associated to this modification.
Thus we have constructed an object
$$
\eta:=(\sta(C/S, (x_i), f),(\xi_i))\quad\in\quad\Mb_{g,n}(S)\times\Ws^n(S)=
\Mt_{g,n}(S)
\quad,
$$
which we define to be the image of $(C/S, (x_i), f)$ by $\stt$.
\epf

Next we recall the definition of some tautological classes
on $\Mb_{g,n}(V,\beta)$ and at the same time introduce notation
which we will use in the sequel.

\bdefin
\label{phipsi}
\barabic
\item
Let $g,n\geq 0$ and let $(X_{g,n}(V,\beta)/\Mb_{g,n}(V,\beta),(x_i),f)$
be the universal object over $\Mb_{g,n}(V,\beta)$. For $0\leq i\leq n-1$
we define the line bundle
$$
L_{i;g,n}(V,\beta):=x_i^*\omega_{X_{g,n}(V,\beta)/\Mb_{g,n}(V,\beta)}\quad.
$$
over $\Mb_{g,n}(V,\beta)$.
\item
If $(g,n)$ is in the stable range and $V$ is a point then
$\Mb_{g,n}(V,\beta)=\Mb_{g,n}$ and we write 
$L_{i;g,n}$ instead of $L_{i;g,n}(V,\beta)$.
\item
For $g,n\geq 0$ we let 
$$
\psi_{i;g,n}(V,\beta):=c_1(L_{i;g,n}(V,\beta))
\quad\in A_*(\Mb_{g,n}(V,\beta))
$$
and for $(g,n)$ in the stable range we let
$$
\phi_{i;g,n}(V,\beta):=\sta^*c_1(L_{i;g,n})
\quad\in A_*(\Mb_{g,n}(V,\beta))
$$
In cases where no confusion is likely to arise
we often simply write $\psi_i$ and $\phi_i$ instead
of $\psi_{i;g,n}(V,\beta)$ and $\phi_{i;g,n}(V,\beta)$.
\earabic
\edefin

Observe that for $(g,n)$ in the stable range and
for general $(V,\beta)$ 
due to the fact that $\sta:\Mb_{g,n}(V,\beta)\to\Mb_{g,n}$
contracts modifications,
the line bundles $\sta^*L_{i;g,n}$
and $L_{i;g,n}(V,\beta)$ are not isomorphic.
Therefore in general the classes $\psi_i$ and $\phi_i$ are 
different and it is not clear whether $\psi_i$ is a pull
back by $\sta$ of any class in $\Mb_{g,n}$. 
One of the advantages of the morphism $\stt$ over
the morphism $\sta$ is that, as we will see in a moment,
for any $(g,n)$ in the half-stable range
there are
classes $\psit_i\in A_*\Mt_{g,n}$
whose pull back via $\stt$ are the classes $\psi_i$.

\bdefin
\label{psit}
\barabic
\item
We define the following elements in $A_*\Moz=R[t_0]$:
$$
\phit_{0;0,2}:=-t_0-t_1
\quad\text{,}\quad
\phit_{1;0,2}:=t_0
$$
and
$$
\psit_{0;0,2}:=-t_0
\quad\text{,}\quad
\psit_{1;0,2}:=t_0+t_1
\quad.
$$
\item
Let $(g,n)$ be in the stable range. 
For $0\leq i\leq n-1$  we define the elements
$$
\phit_{i;g,n}:=c_1(L_{i;g,n})\tensor 1\tensor\dots\tensor 1
$$
and
$$
\psit_{i;g,n}:=\phit_{i;g,n}+t^{(i)}_1
$$
in $A_*\Mt_{g,n}=A_*\Mb_{g,n}\tensor R^{\tensor n}$, where
$$
t^{(i)}_1:=
[\Mb_{g,n}]\tensor 
\underset{\text{$(i-1)$ times}}
{\underbrace{1\tensor\dots\tensor 1}}
\tensor t_1\tensor 
\underset{\text{$(n-i)$ times}}
{\underbrace{1\tensor\dots\tensor 1}}
\quad \in A_*\Mb_{g,n}\tensor R^{\tensor n} \quad.
$$
\earabic
Again we often write $\phit_i$ and $\psit_i$ instead of
$\phit_{i;g,n}$ and $\psit_{i;g,n}$ if this does not give
rise to confusion.
\edefin

It is clear from the definitions that we have 
$\stt^*(\phit_i)=\phi_i$ for $(g,n)$ in the stable
range. The classes $\psit_i$ satisfy an analogous equation:
\bpr
\label{stt*psi}
Let $(g,n)$ be in the half-stable range. Then
we have
$$
\stt^*(\psit_i)=\psi_i
\quad\text{for $0\leq i\leq n-1$.}
$$
\epr

\bpf
For $(g,n)$ in the stable range the  equality  follows from Lemma \ref{psi}.
Assume now $(g,n)=(0,2)$.
Let $\Upsilon=(C/S,(x_0,x_1),f:C\to V)$ 
be an $S$-valued point of $\Mb_{0,2}(V,\beta)$.
From the definition of $\stt$ and
the proof of Proposition \ref{M02} it follows that $\stt$ maps this
$S$-valued point to $\eta=(L,\xi)\in\BGm\times\Ws$ where 
$L=x_0^*\Oo_C(x_0)=x_0^*\omega_{C/S}^\mo$.
Thus we have $\Upsilon^*\stt^*t_0=\eta^*t_0=c_1(L)=-\Upsilon^*\psi_0$.
This shows $\stt^*\psit_0=\psi_0$.
Interchanging the role of $x_0,x_1$  we get 
$\psi_1=\stt^*\iota(\psit_0)=\stt^*\psit_1$
where $\iota$ is the involution defined in Proposition \ref{iota}.
\epf

The following proposition will play an important role in
\S\ref{section gravitation}.

\bpr
\label{important}
\balph
\item
In $A_*\Mt_{0,2}$ we have the following equality:
$$
\psit_0^{d_0}\psit_1^{d_1}=
-\psit_0^{d_0+1}\psit_1^{d_1-1}+
\psit_0^{d_0}\comp_1\psit_1^{d_1-1}
\quad\text{for $d_0\geq 0$, $d_1>0$.}
$$
\item
Let $n\geq 2$, let $d_0,\dots,d_n,e_0,\dots,e_n\geq 0$, 
and assume $d_j>0$ for some $j\in\{1,\dots,n\}$. 
Then the following identity holds in $A_*\Mt_{0,n+1}$
$$
\prod_{i=0}^n\psit_i^{d_i}\phit_i^{e_i}=
\prod_{i=0}^n\psit_i^{d_i-\delta_{i,j}}\phit_i^{e_i+\delta_{i,j}}
 + 
\left(
\phit_j^{e_j}
\prod_{i:i\neq j}\psit_i^{d_i}\phit_i^{e_i}
\right)
\comp_j\psit_{1;0,2}^{d_j-1}
$$
where $\psit_i:=\psit_{i;0,n+1}$ and $\phit_i:=\phit_{i;0,n+1}$.
\ealph
\epr

\bpf

For the proof of (a) notice that we have
$$
\psit_0^{d_0}\psit_1^{d_1}=
(-t_0)^{d_0}(t_0+t_1)^{d_1}=
-(-t_0)^{d_0+1}(t_0+t_1)^{d_1-1}+
(-t_0)^{d_0}t_1(t_0+t_1)^{d_1-1}
$$
The stated equality now follows by Lemma \ref{identity}.
The proof of (b) is analogous to the proof of Corollary \ref{cor psi}.
\epf

\section{Gravitational correlators}
\label{section gravitation}

Let $V$ be a smooth projective variety.
For $g,n\geq 0$ and $\beta\in B(V)$ let 
$\Mb_{g,n}(V,\beta)$ be the moduli stack of stable $n$-pointed maps
of genus $g$ and class $\beta$ into $V$ (cf. beginning of \S\ref{section stt}). 
In \cite{BF}, Behrend and Fantechi have constructed canonical element
$J_{g,n}(V,\beta)\in A_*\Mb_{g,n}(V,\beta)$,
called the ``virtual fundamental class''.

Assume now that $(g,n)$ is in the stable range, i.e. that $n>2(1-g)$.
Then the moduli space $\Mb_{g,n}$ is not empty and
there exist natural morphisms
\begin{align*}
\ev:\Mb_{g,n}(V,\beta) &\to V^n
\\
\sta:\Mb_{g,n}(V,\beta) &\to\Mb_{g,n}
\end{align*}
called the ``evaluation morphism'' and the ``stabilisation morphism''
respectively (cf. \cite{Manin} Ch.V \S4).
The push-forward
$$
I_{g,n}(V,\beta):=(\ev,\sta)_*J_{g,n}(V,\beta)
\quad\in A_*(V^n\times\Mb_{g,n})
$$
of the virtual fundamental class
is called the ``Gromov-Witten correspondence''.
For cohomology classes $\alpha_0,\dots,\alpha_{n-1}\in H^*(V,\Q)$
one defines the Gromov-Witten class
$$
I_{g,n,\beta}(\alpha_0,\dots,\alpha_{n-1}):=
\PD^\mo p_{2*}(p_1^*(\alpha_0\tensor\dots\tensor\alpha_{n-1})\cap I_{g,n}(V,\beta))
\quad\in H^*(\Mb_{g,n})
$$
where $p_1:V^n\times\Mb_{g,n}\to V^n$ and $p_2:V^n\times\Mb_{g,n}\to \Mb_{g,n}$
are the projections, and $\PD: H^*(\Mb_{g,n})\isomto H_*(\Mb_{g,n})$
is the Poincar\'e duality map. 
Let 
$
I^\dual_{0,n,\beta}:
H_*(\Mb_{0,n},\Q)\to H_*(V,\Q)^{\tensor n}
$
be dual to the map 
$
\alpha_0\tensor\dots\tensor\alpha_{n-1}\mapsto
I_{g,n,\beta}(\alpha_0,\dots,\alpha_{n-1})
$
and let $I^\dual_n$ be the composition of
$$
\sum_{\beta\in B(V)}I^\dual_{0,n,\beta}q^\beta:
H_*(\Mb_{0,n},\Q)\to H_*(V,\Lambda)^{\tensor n}
$$
with the natural map $A_*(\Mb_{0,n})\to H_*(\Mb_{0,n},\Q)$
The following theorem is well-known.

\bthm
\label{Idual}
Let $\Lambda$ be the topological
ring defined at the beginning of \S\ref{section stt}.
Let $H$ be the $\Lambda$-module $H_*(V,\Lambda)$,
endowed with the natural $\Z/2$-grading.
Let $b:H\tensor H\to\Lambda$ be the $\Lambda$-linear
continuation of the Poincar\'e pairing.
Then the morphisms $I^\dual(1):=0$ and
$$
I^\dual(n):=I^\dual_{n+1}:A\Mb_0(n)\to\E^s[H,b](n)
\quad\text{for $n\geq 2$}
$$ 
define a morphism $I^\dual:A\Mb_0\to\E^s[H,b]$ of cyclic operads.
\ethm

Recall the definition of Gromov-Witten invariants:
For $\alpha_0,\dots,\alpha_n\in H^*(V,\Q)$ one sets
$$
\langle
\alpha_0,\dots,\alpha_n
\rangle_{0,\beta}:=
\int_{J_{0,n}(V,\beta)}\ev_0^*\alpha_0\cup\dots\cup\ev_n^*\alpha_n
\quad\in\Q
$$
where $\ev_i:\Mb_{0,n+1}(V,\beta)\to V$ is the $i$-th evaluation map.
The following Lemma is easy to check:

\blem
\label{GW-inv}
In terms of the morphism $I^\dual:A\Mb_0\to\E^s[H,b]$
the Gromov-Witten invariants can be expressed as follows:
$$
\sum_{\beta\in B(V)}
\langle
\alpha_0,\dots,\alpha_n
\rangle_{0,\beta}\ q^\beta=
\Bigl\langle
I^\dual(n)([\Mb_{0,n+1}])\ ,\
\alpha_0\tensor\dots\tensor\alpha_n
\Bigr\rangle
$$
where the brackets on the right hand side 
are induced by the pairing $H_*(V,\Q)\tensor H^*(V,\Q)\to\Q$. 
\elem

Recall (cf. \S\ref{section stt}) that for all $n>1-g$
there is a morphism
$$
\stt:\Mb_{g,n}(V,\beta)\to \Mt_{g,n}:=\Mb_{g,n}\times\Ws^n
$$
which is analogous to the stabilisation  morphism $\sta$.
By means of this morphism, we would like to construct 
an analogue of the morphism $I^\dual$. However it is not
clear to us how to do this, since the morphisms $\stt$ are
not proper. Still it seems reasonable to make the following
conjecture:

\bcon
\label{conj1}
There is a natural morphism $\It^\dual:A\Mt_0\to\E^s[H,b]$
of cyclic operads, which is analogous to the morphism $I^\dual$
and has the following properties:
\balph
\item
For each $n\geq 1$ the image of $\It^\dual(n)$ lies in the even
part of $H^{\tensor(n+1)}$.
\item
The morphism $I^\dual$ from Theorem \ref{Idual} factorises through
$\It^\dual$. More precisely we have $I^\dual=\It^\dual\comp F$
where $F:A\Mb_0\to A\Mt_0$ is the morphism from Example \ref{AMb to AMt}.
\ealph
\econ

We expect that  also an analogue of Lemma \ref{GW-inv} holds, where
Gromov-Witten invariants are replaced by gravitational correlators
(also called Gromov-Witten invariants with gravitational descendents)
and $I^\dual$ is replaced by $\It^\dual$.
Before formulating a conjecture, let us recall the definition of
(generalised) gravitational correlators. 
Let $\phi_i,\psi_i\in A^*(\Mb_{g,n+1}(V,\beta))$ 
be the classes defined in Definition \ref{phipsi}.
To cohomology classes $\alpha_0,\dots,\alpha_n\in H^*(V,\Q)$
and integers $d_i,e_i\geq 0$ one associates the numbers
$$
\langle
\tau_{d_0,e_0}\alpha_0,\dots,\tau_{d_n,e_n}\alpha_n
\rangle_{g,\beta}:=
\int_{J_{g,n+1}(V,\beta)}
\psi_0^{d_0}\phi_0^{e_0}\cup\ev_0^*\alpha_0
\cup\dots\cup
\psi_n^{d_n}\phi_n^{e_n}\cup\ev_n^*\alpha_n
\quad\in\Q
$$
called {\em generalised gravitational correlators}. 
If all $e_i$ vanish, the numbers are  called {\em gravitational
correlators}.
It is customary to write $\tau_d$ instead of $\tau_{d,0}$.
Notice that the $\phi_i$ are defined only for $(g,n)$ in the stable range
since otherwise $\Mb_{g,n}$ is empty.
In contrast, $\psi_i$ is defined for any $g,n\geq 0$.
So while the above definition of 
generalised gravitational correlators
makes sense only in the stable range, 
gravitational correlators proper exist for all $g,n$.

\bcon
\label{conj2}
Let $\phit_i,\psit_i\in A_*\Mt_{g,n+1}$ be defined as in Definition
\ref{psit}.
The morphism $\It^\dual$ from Conjecture \ref{conj1}
has the property that
$$
\sum_{\beta\in B(V)}
\langle
\tau_{d_0,e_0}\alpha_0,\dots,\tau_{d_n,e_n}\alpha_n
\rangle_{0,\beta}\ q^\beta=
\Bigl\langle
\It^\dual(n)\left(\prod_{i=0}^n\psit_i^{d_i}\phit_i^{e_i}\right)\ ,\
\alpha_0\tensor\dots\tensor\alpha_n
\Bigr\rangle
$$
for all $n\geq 1$, $\alpha_0,\dots,\alpha_n\in H^*(V,\Q)$. 
\econ

In the remaining of this paragraph we assume Conjectures \ref{conj1}
and \ref{conj2}  to be true.
As a first application we will show how the following 
theorem of Kontsevich and Manin
(cf. \cite{KM} Thm. 1.2, \cite{Manin} VI. Thm. 6.2)
can be derived from the conjectures.

\bthm
\label{thm KM}
Let $n>1$, let $d_0,\dots,d_n\geq 0$, $e_0,\dots,e_n\geq 0$,
$\beta\in B(V)$ and let $\alpha_0,\dots,\alpha_n\in H^*(V,\Q)$
be homogeneous elements.
Let $(\Delta_a)$ be a
basis of the $\Q$-vector space $H^*(V,\Q)$
and let $(\Delta^a)$ be its Poincar\'e dual basis.
Assume that $d_j>0$ for some $j$.
Then we have
\begin{multline*}
\Bigl\langle 
\tau_{d_i,e_i}\alpha_i
\st i=0,\dots,n
\Bigr\rangle_{0,\beta}
=
\Bigl\langle 
\tau_{d_i-\delta_{i,j},e_i+\delta_{i,j}}\alpha_i
\st i=0,\dots,n
\Bigr\rangle_{0,\beta}
+
\\
+
\sum_{a,\beta_1+\beta_2=\beta}
(-1)^{P} 
\Bigl\langle 
\tau_{d_j-1}\alpha_j,
\tau_0\Delta^a
\Bigr\rangle_{0,\beta_1}
\Bigl\langle 
\tau_{0,e_j}\Delta_a,
\tau_{d_i,e_i}\alpha_i
\st i=0,\dots,j-1,j+1,\dots,n
\Bigr\rangle_{0,\beta_2}
\end{multline*}
where $P=\sum_{i=0}^{j-1}|\alpha_i||\alpha_j|$.
\ethm

Theorem \ref{thm KM} can be deduced from
Conjectures \ref{conj1} and \ref{conj2} as follows.
Assume that these conjectures are true.
For convenience we write 
$$
\langle \eta;\alpha_1,\dots,\alpha_n \rangle
:=
\langle
\It^\dual(n)(\eta),\alpha_0\tensor\dots\tensor\alpha_n
\rangle
\quad.
$$
Since the gravitational correlators are symmetric (in the sense of super mathematics)
in its entries, the case $j=0$ follows from the case $j=1$. Thus
we may assume $j>0$.
By Proposition \ref{important} we have 
$$
\prod_{i=0}^n\psit_i^{d_i}\phit_i^{e_i}=
\prod_{i=0}^n\psit_i^{d_i-\delta_{i,j}}\phit_i^{e_i+\delta_{i,j}}
 + 
\left(
\phit_j^{e_j}
\prod_{i:i\neq j}\psit_i^{d_i}\phit_i^{e_i}
\right)
\comp_j\psit_{1;0,2}^{d_j-1}
\quad.
$$
Without loss of generality we may assume that the $\Delta_\nu$, $\Delta^\nu$
are homogeneous.
Lemma \ref{vowa} from the Appendix implies
\begin{multline*}
\left\langle
\left(
\phit_j^{e_j}
\prod_{i:i\neq j}\psit_i^{d_i}\phit_i^{e_i}
\right)
\comp_j\psit_{1;0,2}^{d_j-1}
;
\alpha_0,\dots,\alpha_n
\right\rangle
=\\=
\sum_\nu
(-1)^{N_\nu}
\left\langle
\left(
\phit_j^{e_j}
\prod_{i:i\neq j}\psit_i^{d_i}\phit_i^{e_i}
\right)
;
\alpha_0,\dots,\alpha_{j-1},
\Delta_\nu,
\alpha_{j+1},\dots,\alpha_n
\right\rangle
\left\langle
\psit_{1;0,2}^{d_j-1}
;
\Delta^\nu,\alpha_j
\right\rangle
\end{multline*}
where 
$$
N_\nu=
P+
|\alpha_j||\Delta^\nu|+
\sum_{i=0}^{j-1}|\alpha_i||\Delta_\nu|
$$
By means of Conjecture \ref{conj2} 
the Theorem follows immediately.

Up to this point we have used Proposition \ref{important} only
for $n>2$, i.e. in the stable range. Let us see what it looks
like for $n=2$. Recall that in this case the formula reads
$$
\psit_0^{d_0}\psit_1^{d_1}=
-\psit_0^{d_0+1}\psit_1^{d_1-1}+
\psit_0^{d_0}\comp_1\psit_1^{d_1-1}
\quad\text{for $d_0\geq 0$, $d_1>0$.}
$$
Applying Lemma 
\ref{vowa} from the Appendix and assuming
Conjectures \ref{conj1}, \ref{conj2} we get

\bpr
\label{formula}
Let $d_0\geq 0$, $d_1> 0$,
$\beta\in B(V)$ and let $\alpha_1,\alpha_2\in H^*(V,\Q)$
be homogeneous elements.
Let $(\Delta_a)$ be a
basis of the $\Q$-vector space $H^*(V,\Q)$
and let $(\Delta^a)$ be its Poincar\'e dual basis
(i.e. $\Delta_a\cup\Delta^b=\delta_{a,b}$).
Then we have
\begin{multline*}
\Bigl\langle 
\tau_{d_0}\alpha_0,
\tau_{d_1}\alpha_1
\Bigr\rangle_{0,\beta}
=
-
\Bigl\langle 
\tau_{d_0+1}\alpha_0,
\tau_{d_1-1}\alpha_1
\Bigr\rangle_{0,\beta}
+
\\
+
\sum_{a,\beta_1+\beta_2=\beta}
(-1)^{|\alpha_0||\alpha_1|}
\Bigl\langle 
\tau_{d_0}\alpha_0,
\tau_{0}\Delta_a
\Bigr\rangle_{0,\beta_1}
\Bigl\langle 
\tau_0\Delta^a,
\tau_{d_1-1}\alpha_1
\Bigr\rangle_{0,\beta_2}
\quad.
\end{multline*}
\epr

The formula in Proposition \ref{formula} 
can be proven by other means
then via Conjectures \ref{conj1}, \ref{conj2}. 
In fact A. Givental assured me that it is well-known
since the very inception of the theory (Dijkgraaf-Witten and Dubrovin, 1991-1992).
For example it is implied by the formula (for $t=0$)
$$
V_{\alpha,\beta}(t,q,x,y)
=
\frac{1}{x+y}
\sum_{\epsilon,\epsilon'}
\eta^{\epsilon,\epsilon'}
S_{\epsilon,\alpha}(t,q,x)
S_{\epsilon',\beta}(t,q,y)
$$
on the bottom of the 13-th page of
Givental's paper \cite{Givental}.

The following proposition reduces the amount of information contained
in the morphism $\It^\dual:A\Mt_0\to\E^s[H,b]$.

\bpr
\label{final}
Assume that Conjectures \ref{conj1}, \ref{conj2} are true.
Then the morphism $\It^\dual:A\Mt_0\to\E^s[H,b]$ of cyclic operads
is completely determined by all two-point gravitational correlators
thogether with all elements
$$
\left\langle
\xi\prod_{i=0}^n\psit_i^{k_i};\alpha_0,\dots,\alpha_n
\right\rangle
\quad \in\Lambda
$$
where $n\geq 2$, $k_i\geq 0$, $\xi\in A_*\Mb_{0,n+1}$,
$\alpha_i\in H^*(V,\Q)$.
\epr

\bpf
Clearly, the knowledge of
$\It^\dual:A\Mt_0\to\E^s[H,b]$
amounts to the same as the knowledge of the elements
$$
\langle
\eta; \alpha_0,\dots,\alpha_n
\rangle
\quad\in\Lambda
$$
for all $n\geq 1$, $\eta\in A\Mt_0(n)$, $\alpha_i\in H^*(V,\Q)$.

Assume first that $n=1$. By Corollary \ref{basis} the element $\eta$
is a linear combination of elements of the form
$$
\psit_1^{a_1}\odot\dots\odot\psit_1^{a_r}=
(\dots(\psit_1^{a_1}\comp_1\psit_1^{a_2})\comp_1\dots)\comp_1\psit_1^{a_r}
$$
where $r\geq 1$, $a_i\geq 0$ and $\psit_1=t_0+t_1$.
Applying Lemma \ref{vowa} of the Appendix, it follows that the element 
$$
\langle
(\dots(\psit_1^{a_1}\comp_1\psit_1^{a_2})\comp_1\dots)\comp_1\psit_1^{a_r};
\alpha_0,\dots,\alpha_n
\rangle
\quad\in\Lambda
$$
can be written as an expression involving only gravitational two-point
correlators. This proves the proposition in the case $n=1$.

Now let $n\geq 2$.
Any element of $A\Mt_0(n)$ is a linear combination of elements of the
form
$$
\eta=\xi\tensor w_0\tensor\dots\tensor w_n
$$
where $\xi\in A\Mb_0(n)$ and $w_i\in R$.
With the notation on page \pageref{q'} we can write
$
\eta=\xi\prod_{i=0}^nq'_i([\Mb_{0,n+1}]\tensor w_i)
\quad.
$
Since $w_i$ does not depend on $t_0$, we have 
$[\Mb_{0,n+1}]\tensor w_i=w_i|_{t_0=\phi_i}$ in the notation of
Proposition \ref{cK*}, and therefore we can
write
$$
\eta=\xi\prod_{i=0}^nq'_i(w_i|_{t_0=\phi_i})
\quad.
$$
By Corollary \ref{basis} every $w_i$ is a linear combination of
elements of the form $\psit_1^{a_1}\odot\dots\odot\psit_1^{a_r}$,
and it is easy to check that for all $x_0,\dots,x_n,y\in R[t_0]$
we have
$$
\xi\left(\prod_{i:i\neq j}q'_i(x_i|_{t_0=\phi_i})\right)
q'_j((x_j\odot y)|_{t_0=\phi_j})=
\left(\xi\prod_{i:i\neq j}q'_i(x_i|_{t_0=\phi_i})\right)\comp_j y
\quad.
$$
This shows that with the help of Lemma \ref{vowa} from the Appendix
we can again write the element $\langle \eta; \alpha_0,\dots,\alpha_n\rangle$
as an expression involving only gravitational two-point correlators
and elements of the form 
$\langle \gamma; \alpha_0,\dots,\alpha_n\rangle$
where
$$
\gamma=\xi\prod_{i=0}^nq'(\psit_1^{k_i}|_{t_0=\phi_i})=
\xi\prod_{i=0}^n\psit_i^{k_i}
$$
for some $k_i\geq 0$.
\epf

\brem
\label{finalrem}
With the notation of Proposition \ref{final}
it seems reasonable to conjecture that for $n\geq 2$ we have
$$
\left\langle
\xi\prod_{i=0}^n\psit_i^{k_i};\alpha_0,\dots,\alpha_n
\right\rangle
=
\sum_{\beta\in B(V)}q^\beta
\int_\xi I_{0,n+1}(V,\beta;k_0,\dots,k_n)(\alpha_0\tensor\dots\tensor\alpha_n)
$$
where the maps 
$I_{g,n+1}(V,\beta;k_0,\dots,k_n):H^*(V^{n+1},\Q)\to H^*(\Mb_{g,n+1},\Q)$
are the so-called {\em cohomological correlators} defined in
\cite{Manin}, Ch. VI \S2.2 (2.4).
Granting this, Proposition \ref{final}
says that the $A\Mt_0$-algebra structure
on $H^*(V,\Lambda)$ amounts to the same thing as full gravitational quantum
cohomology of $V$.
\erem

\section{Appendix}
\label{appendix}

We place ourselves in the situation of Example \ref{EsMb}.
Thus let $\kb$ be a ring, let $M$ be a $\Z/2$-graded
$\kb$-module and let $b:M\times M\to \kb$ be an even bilinear
form on $M$.

The bilinear form $b$ induces an even bilinear form $\langle\ ,\ \rangle$  on 
$M^{\tensor n}$ which
is characterised as follows:
$$
\langle
v_1\tensor\dots\tensor v_n,
\alpha_1\tensor\dots\tensor \alpha_n
\rangle
=
(-1)^S\prod_{i=1}^nb(v_i,\alpha_i)
$$
for homogeneous $v_i,\alpha_i\in M$, where
$
S=\sum_{i>j}|v_i||\alpha_j|
$. In other words, $(-1)^S$ is the sign which arises from the sign
rule when one reorders the elements $v_1,\dots,v_n,\alpha_1,\dots,\alpha_n$
to $v_1,\alpha_1,\dots,v_n,\alpha_n$.

\blem
\label{vowa}
In the above situation assume that $M$ is
a finitely generated free $\kb$-module and that there are
basis $(\Delta_\nu)$, $(\Delta^\nu)$ of $M$ consisting
of homogeneous elements, such that 
$b(\Delta_\nu,\Delta^{\nu'})=\delta_{\nu,\nu'}$.
Let $m,n\geq 1$, $1\leq j\leq m$, and let $v\in M^{\tensor (m+1)}$, 
$w\in M^{\tensor (n+1)}$ be even elements.
Then for any
$\alpha=\alpha_0\tensor\dots\tensor\alpha_{m+n-1}\in M^{\tensor (m+n)}$
the following equation holds:
$$
\langle v\comp_j w, \alpha \rangle=
\sum_\nu
(-1)^{N_\nu}
\langle
v,
\alpha_0\tensor\dots\tensor\alpha_{j-1}
\tensor\Delta_\nu\tensor
\alpha_{j+n}\tensor\dots\tensor\alpha_{m+n-1}
\rangle
\langle
w,
\Delta^\nu\tensor 
\alpha_{j}\tensor\dots\tensor\alpha_{j+n-1}
\rangle
$$
where
$\comp_j$ is the composition morphism in the operad 
$\E^s(M,b)$ defined in Example \ref{EsMb} and where
$$
N_\nu:=
\sum_{i=0}^{j-1}|\alpha_i||\Delta_\nu|+
\sum_{i=j}^{j+n-1}|\alpha_i||\Delta^\nu|+
\left(
\sum_{p=j}^{j+n-1}|\alpha_p|
\right)
\left(
\sum_{q=0}^{j-1}|\alpha_q|
\right)
\equiv |\Delta_\nu|
\quad\text{mod $(2)$}
\quad.
$$
\elem

\bpf
We may assume that
$v=v_0\tensor\dots\tensor v_{m}$
and 
$w=w_0\tensor\dots\tensor w_{n}$
for homogeneous $v_i$, $w_i$.
By definition we have
$$
v\comp_j w=
\pm
b(v_j, w_0)
v_0\tensor\dots\tensor v_{j-1}\tensor
w_1\tensor\dots\tensor w_n\tensor
v_{j+1}\tensor\dots\tensor v_m
$$
where the sign arises from interchanging $v_{j+1}\tensor\dots\tensor v_m$
and $w_1\tensor\dots\tensor w_n$.
Therefore 
\begin{multline*}
\langle v\comp_j w, \alpha \rangle
=
\epsilon_1
b(v_j,w_0)
\langle
v_0\tensor\dots\tensor v_{j-1}\tensor v_{j+1}\tensor\dots\tensor v_m
\ ,\
\alpha_0\tensor\dots\tensor \alpha_{j-1}\tensor
\alpha_{j+n}\tensor\dots\tensor \alpha_{m+n-1}
\rangle
\\
\langle
w_1\tensor\dots\tensor w_n
\ ,\
\alpha_{j}\tensor\dots\tensor \alpha_{j+n-1}
\rangle
\end{multline*}
where the sign $\epsilon_1\in\{\pm 1\}$ comes from interchanging 
$w_0$ with $v_{j+1},\dots,v_m$, 
interchanging
$w_1,\dots,w_n$ with $\alpha_0,\dots,\alpha_{j-1}$
and 
interchanging
$w_1,\dots,w_n,\alpha_j,\dots,\alpha_{j+n-1}$ with the elements
$\alpha_{j+n},\dots,\alpha_{m+n-1}$.
Observe that we have
$$
b(v_j,w_0)=
\sum_\nu b(v_j,\Delta_\nu)b(w_0,\Delta^\nu)
$$
and consequently
\begin{multline*}
\langle v\comp_j w, \alpha \rangle
=\\=
\sum_\nu
\epsilon(\nu)
\langle
v,
\alpha_0\tensor\dots\tensor\alpha_{j-1}
\tensor\Delta_\nu\tensor
\alpha_{j+n}\tensor\dots\tensor\alpha_{m+n-1}
\rangle
\langle
w,
\Delta^\nu\tensor 
\alpha_{j}\tensor\dots\tensor\alpha_{j+n-1}
\rangle
\end{multline*}
where $\epsilon(\nu)=\epsilon_1\epsilon_2(\nu)$ and $\epsilon_2(\nu)$
is the sign that comes from
moving $\Delta^\nu$ pass the elements
$v_{j+1},\dots,v_m,\alpha_0,\dots,\alpha_{j-1}$,
and moving $\Delta_\nu$ pass $w_1,\dots,w_n$.

To simplify $\epsilon(\nu)$ for given $\nu$, we may
assume that $\Delta^\nu$, $\Delta_\nu$, $v_j$, $w_0$ all have the
same parity, since otherwise the corresponding summand is zero
anyway. Furthermore, since $w$ is assumed to be even,
interchanging $w_1,\dots,w_n$ with something
leads to the same sign 
as interchanging $w_0$ with the same thing.

Therefore the sign $\epsilon(\nu)$ is the same
as the one which arises from simply interchanging
$\Delta_\nu$ with itself and interchanging 
$\Delta^\nu,\alpha_j,\dots,\alpha_{j+n-1}$ with $\alpha_{j+n},\dots,\alpha_{m+n-1}$.
The last interchange has no effect since for given $\nu$ we may assume that
the parity of $\Delta^\nu\tensor\alpha_j\tensor\dots\tensor\alpha_{j+n-1}$
is that of the even element $w$, 
since otherwise the corresponding summand vanishes.
So finally we have $\epsilon(\nu)=(-1)^{|\Delta_\nu|}$.

It remains to show that $N_\nu\equiv|\Delta_\nu|$ mod $(2)$.
But this follows immediately from
$$
|\Delta^\nu|+\sum_{i=j}^{j+n-1}|\alpha_i|\equiv |w|\equiv 0
\quad\text{mod $2$.}
$$
\epf


\end{document}